\magnification=\magstep1
\input epsf

~~\vskip .5in

\vskip .5in

\centerline{\bf The Structure and Singularities of Arc Complexes}
\vskip .2in

\centerline {\bf R. C. Penner}

\centerline{Departments of Mathematics and Physics/Astronomy}

\centerline{University of Southern California}

\centerline {Los Angeles, CA 90089}

\vskip .2in

\centerline{}

\vskip .2in

\leftskip .8in\rightskip .8in

\noindent {\bf Abstract}~~A classical combinatorial fact is 
that the
simplicial complex consisting of disjointly embedded chords in a convex 
planar polygon is a sphere.  For any
surface $F$ with non-empty boundary, there is an analogous complex $Arc(F)$ 
consisting of suitable
equivalence classes of arcs in $F$ connecting its boundary components.
The main result of this paper is 
the determination
of those arc complexes $Arc(F)$ that are also spherical.  
This classification has consequences for
Riemann's moduli space of curves via its
known identification with an analogous arc complex
in the punctured case with no boundary.  Namely, the essential singularities
of the natural cellular compactification of Riemann's moduli space can be described.

\leftskip=0ex\rightskip=0ex

~~\vskip .2in

\centerline{ \bf Introduction}

\vskip .2in

\noindent The projectivized lamination space of a surface $S$
introduced by Thurston [27] in the 1970's
is a fundamental tool in dynamics, topology, and geometry
in low dimensions,
and it contains two characteristic
subspaces which have also played roles in subsequent
developments: The ``curve complex'' of $S$
is the subspace of projective laminations
of compact support in $S$ where each leaf is a simple closed curve, 
and the ``pre-arc
complex'' $Arc'(S)$
of a punctured surface $S$ without boundary
is the subspace of projective laminations
where each leaf is an arc asymptotic in both directions to the punctures.
                                                                                
\vskip .1in
                                                                                
\noindent
Harvey [9] introduced the
curve complex to give a real-analytic manifold with boundary which
compactifies the Teichm\"uller space so that the mapping class group acts
discontinuously on the compactification.  
Studying the connectivity of the curve complex, he
then used this to prove that the mapping class groups are
virtual duality groups and to provide
a lower bound to their virtual cohomological
dimensions.
                                                                                
\vskip .1in
                                                                                
\noindent
Whereas the
mapping class group action on the space of projective laminations 
has a non-Hausdorff quotient, the action on $Arc'(S)$
is discontinuous (as is the action on the curve complex), 
and the quotient $Arc(S)$ is the ``arc complex'' of the
underlying punctured surface $S$.  $Arc'(S)$ contains a further subspace
defined by the condition that any essential closed curve
in $S$ has non-zero transverse measure, and the
quotient is denoted $Arc_\#(S)\subseteq Arc(S)$.

\vskip .1in
                                                                                
\noindent
Mumford realized that deep analytic results of Strebel [26] and Hubbard-Masur
[10] imply that $Arc_\#(S)$ is real-analytically homeomorphic to the
product of Riemann's moduli space of $S$ with an open simplex of dimension
$s-1$, where $s\geq 1$ is the number of punctures of $S$.  Two related
treatments in [17] and subsequently by Bowditch-Epstein in
[1] also gave the same result in the setting of hyperbolic
rather than conformal geometry.  Immediate consequences include
a presentation for the fundamental path groupoid of Riemann's moduli
space [17] and a method of integration [23]
over moduli space.
Furthermore, it is not difficult to see that
$Arc(S)$ is compact, so it forms a natural compactification of the space
$Arc_\#(S)$, which is identified with the product of Riemann's moduli
space and a simplex.   (See [1,19] and $\S$6 for more details about
this and related compactifications.)
                                                                                
\vskip .1in

\noindent Harer [6] presented Mumford's argument and
used this identification together with coordinates on $Arc'(S)$
from [18] (see Appendix B) to give an equivariant deformation
retraction of the Teichm\"uller space of $S$ to a lower-dimensional spine
and to provide a Borel-Serre type
bordification. The boundary of this bordification was
shown to be homotopy equivalent to a wedge of spheres, giving another
proof that mapping class groups are virtual duality groups.  In fact,
the dimension of the boundary agreed with Harvey's earlier bound,
thus providing an exact calculation of the virtual cohomological dimension.
In [7], Harer further studied the action of the mapping
class group on various highly connected curve and pre-arc
complexes using
spectral sequence techniques relating the desired homology with that
of stabilizers
to prove that the $k$th homology of the mapping class group is independent
of the genus when it
is large compared to $k$.   This stable range for homology was subsequently
sharpened by Ivanov [12,13] using analogous techniques.
                                                                                
\vskip .1in
                                                                                
\noindent There have been many further applications of arc complexes including
the calculation [8,22] of the virtual Euler characteristic of Riemann's
moduli space.  In particular, the latter treatment uses the
natural identification induced by Poincar\'e duality in $S$
of $Arc_\#(S)$ with a complex of graphs embedded as spines of $S$ together
with matrix model techniques from high energy physics, which were partly
rediscovered in [8].  These graph complex techniques have been further
developed for instance by Kontsevich [15] 
in his proof of Witten's conjecture
relating intersection numbers of tautological classes to an
infinite integrable hierarchy of Korteweg-de Vries equations and 
by Igusa [11] and Mondello [16] who independently prove
that the tautological classes can be described
combinatorially.
                                                                                
\vskip .1in

\noindent
The current paper is primarily concerned with the analogous arc complex
for any {\it bordered} surface $F$ of genus $g\geq 0$ with $s\geq 0$ punctures
and $r\geq 1$ boundary components, where we choose
for each $i=1,\dots ,r$ some number $\delta _i\geq 1$ of
distinguished points on the $i^{\rm th}$
boundary component of $F$ at which arcs may terminate.
One considers equivalence classes of families of
disjoint arcs in the surface
connecting the distinguished points on the boundary in order to define a
corresponding ``arc complex'' $Arc(F)$. (See
$\S$1 for the precise definition.)  In analogy to the punctured
case, there is a subspace of $Arc(F)$
which is closely related
to Riemann's moduli space of the bordered surface (see [20]
for details), but equally to the point, arc complexes of bordered
surfaces, or finite quotients of them, arise as links of simplices
in the natural PL structure of arc complexes of punctured surfaces.

\vskip .1in

\noindent Thus, the global topology of arc complexes for bordered surfaces
studied here
governs the local structure of the arc complex for a punctured surface, that
is, the local structure of a compactification of 
(undecorated) Riemann's moduli space.  (See $\S$6 for the precise statement.)

\vskip .1in
                                                                                
\noindent
The simplest case of a bordered surface is a convex $n$-gon in the plane
(i.e., $g=s=r-1=0$ and $\delta _1=n$), for $n\geq 4$, where the arc complex
is a simplicial complex whose $p$-simplices are sets of $p+1$ disjoint
chords of the polygon and where the chords are required to have their endpoints
at the vertices of the polygon.  We shall give an elementary proof (in $\S$3)
of following well-known result (which Giancarlo Rota once told us was known to Hassler Whitney).
                                                                                
\vskip .2in
                                                                                
\noindent {\bf (Classical) Fact} \it The arc complex of an $n$-gon
is PL-homeomorphic to the sphere of dimension $n-4$.\rm
                                                                                
\vskip .2in
                                                                                
\noindent Given two bordered surfaces $M_1,M_2$, we may consider inclusions $M_1\subseteq M_2$, where the distinguished points and punctures
of $M_1$ map to those of $M_2$, and $M_1$ is a complementary component to an arc
family (possibly an empty arc family) 
in $M_2$.  

\vskip .1in

\noindent Define the {\it type 1 surfaces} to
be the following bordered surfaces: the torus-minus-two-disks (i.e., $g=1,\delta _1+\delta _2=r=2,s=0$), $s$-times-punctured
sphere-minus-$r$-disks with
$r+s=4$ (i.e.,
$g=0$,
$r+s=4$, $\delta _1+\cdots \delta _r=r\geq 1$).  
                                                                                
\vskip .1in
                                                                                
\noindent Here is the main result of this paper:
                                                                                
\vskip .25in

\noindent {\bf Theorem 1}\it ~~The arc complex $Arc(F)$ of a
bordered surface $F$ is PL-homeomorphic to the sphere of dimension
$6g-7+3r+2s+\delta _1+\delta _2+\cdots +\delta _r$ if and only if 
$M\not\subseteq F$ for any type 1 surface $M$.  
In other words, $Arc(F)$ is spherical only
in the following cases:  polygons ($g=s=0$, $r=1$), multiply punctured polygons
($g=0$, $r=1$),
``generalized'' pairs of pants ($g=0$, $r+s=3$), the torus-minus-a-disk ($g=r=1$, $s=0$), and the once-punctured torus-minus-a-disk
($g=r=s=1$).  Only the type 1 surfaces have an arc complex which is 
a PL-manifold other than a sphere. \rm

\vskip .25in

\noindent Sphericity of the arc complex of a polygon is the Classical
Fact, and sphericity for a multiply
punctured polygon was first proven in [19].  (An 
alternate argument for punctured polygons
is given here in $\S$4.)  The sphericity 
question was raised in [19] for any
bordered surface,
and this question is resolved by Theorem~1.  

\vskip .2in

\noindent {\bf Corollary}~~\it The first non-spherical arc complexes
occur in dimension five, and the first non-manifold arc complexes
occur in dimension six.\rm

\vskip .2in

\noindent
In effect, arc complexes are
indeed spherical until some spurious homology arises, which
first occurs in dimension
five.

\vskip .1in

\noindent This paper is organized as follows.  $\S$1 gives the definition of the arc complexes and includes elementary and basic
observations about them as well as the simplest example of a non spherical arc complex for the twice-punctured annulus.  $\S$2 studies
properties of the geometric realizations of partially ordered sets and contains the basic result that $Arc(F)$ is a manifold if the arc complex
of every proper sub-surface of $F$ is spherical.  In $\S$3, we prove that if $F'$ arises from $F$ by adding one more distinguished point on the
boundary, then under suitable hypotheses $Arc(F')$ is PL-homeomorphic to the suspension of $Arc(F)$, and this result is used to
prove the Classical Fact by induction.  In $\S$4, we prove that if $F'$ arises from $F$ by adding one more puncture, then $Arc(F')$ is
spherical provided that $Arc(F)$ is spherical and $F$ has only one 
boundary component; this is used to prove sphericity of arc complexes of
multiply punctured  polygons and the once-punctured
torus-minus-a-disk.  The remaining
counter-examples to sphericity of arc complexes are presented in $\S$5.  Finally, the proof of Theorem~1 is completed in $\S$6, which also
contains a discussion of the general stratified structure of arc complexes as well as Riemann's moduli space and applications.

\vskip .1in

\noindent We explicitly verify the sphericity
of all arc complexes of dimension at most four in Appendix~A.  
Sphericity of arc complexes for generalized
pairs of pants is proven in Appendix~B using techniques [18,24]
from Thurston theory.  The reader could well read these appendices after
$\S$1 to get the flavor of part of the theorem.

\vskip .2in

\noindent {\bf Acknowledgement}~~It is a pleasure to thank Dennis Sullivan for contributing the proof of 
non-sphericity of arc complexes in Proposition~4, for his guidance on this project, and for countless 
exhilarating discussions and late-night phone calls.
I am likewise happy to thank Bob Edwards, Yair Minsky, and Shmuel Weinberger
for helpul remarks.

\vskip .5in

\noindent{\bf 1. Definition of the arc complex and basic properties}

\vskip .2in

\noindent Fix an oriented and smooth bordered surface
of genus $g\geq 0$ with $s\geq 0$ punctures and $r\geq 1$ boundary
components $\partial _1,\ldots ,\partial _r$, and choose in the $i^{\rm th}$ boundary component some number $\delta
_i\geq 1$ of distinguished points $D_i\neq \emptyset$, for
$i=1,2,\dots r$. 

\vskip .1in

\noindent Let $F=F_{g,\vec\delta }^s$ denote the fixed
surface with this data, where $\vec\delta =(\delta _1,\delta _2,\dots \delta _r)$.  We shall
employ the following notation throughout this paper:
$$D=D(F)=D_1\cup\cdots\cup D_r$$ denotes the set of all the distinguished points in the boundary
of $F$, 
$$\Delta =\Delta (F)=\delta _1+\delta _2+\cdots +\delta _r,$$ is the cardinality of $D$, and 
$$N=N(F)=6g-7+3r+2s+\Delta$$
is the dimension of the arc complex.
Furthermore, in the special case that $\Delta =r$ so $\vec\delta$ is an $r$-dimensional vector each of whose
entries is unity, then we shall also sometimes write $F_{g,r}^s=F_{g,\vec\delta}^s$.  In particular, $N\geq 0$ except for the
triangle $F_{0,3}^0$ and once-punctured monogon $F_{0,1}^1$, whose arc complexes are empty.

\vskip .1in

\noindent Define an {\it (essential) arc} in $F$ to be a smooth path $a$
embedded in $F$ whose endpoints lie in $D$ and which meets $\partial F$ transversely, where
we demand that
$a$ is not isotopic rel endpoints to a path lying in $\partial F -D$.  Two arcs
are said to be {\it parallel} if there is an isotopy between them
which fixes $D$ pointwise.  
An {\it arc family} in
$F$ is the isotopy class of a collection of disjointly embedded
essential arcs in
$F$, no two of which are parallel.

\vskip .1in

\noindent Notice that a non-separating arc in a surface is always uniquely determined by its
endpoints.  Furthermore, any separating arc is uniquely determined by its endpoints and the topological type
of its complementary regions by the classification of topological surfaces.  In
particular, any arc in a planar surface whose endpoints coincide is necessarily
separating.
                                                                                
\vskip .1in

\noindent Let us inductively build a simplicial complex $Arc' (F)$ as in the
Classical Fact, where there is one $p$-simplex $\sigma(\alpha )$ for each arc family $\alpha $ in
$F$ of cardinality $p+1$.  The simplicial structure of $\sigma (\alpha )$ is the natural one, where
faces of 
$\sigma (\alpha )$ correspond to sub arc families of $\alpha $.  We begin with a vertex in $Arc'(F)$
for each isotopy class of essential arc in $F$ to define the 0-skeleton.  Having thus inductively
constructed the $(p-1)$-skeleton of $Arc'(F)$, for $p\geq	 1$, let us adjoin a $p$-simplex for each
arc family $\alpha $ consisting of $(p+1)$ essential arcs, where we identify the proper faces of
$\sigma (\alpha )$ with simplices in the $(p-1)$-skeleton in the natural way.  This completes the
inductive definition of $Arc'(F)$.

\vskip .1in

\noindent It is convenient in the sequel to think of this rather more explicitly, as follows.
If $\alpha $ is an arc family, then a {\it weighting} on $\alpha $ is the assignment of 
a non-negative real number, the {\it weight}, to each arc comprising $\alpha $, and the
weighting is {\it positive} if each weight is positive.  A {\it projective weighting} is the
projective class of a weighting on an arc family.  Since the points of a simplex are described by 
``projective non-negative weightings on its vertices'', we
may identify a point of
$\sigma (\alpha )$ with a projective weighting on $\alpha $ in the natural way, so that the
projective positive weightings correspond to the interior $Int~{\sigma }(\alpha )$ of $\sigma 
(\alpha )$.
$Arc'(F)$ itself is thus identified with the collection of all projective positively weighted arc families
in $F$.

\vskip .1in

\noindent Suppose that $\alpha =\{ a_0,\ldots ,\alpha _p\}$ is an arc family.  A useful realization of a positive weighting
$w=(w_i)_0^p$ on the respective components of $\alpha$ is as follows.  One imagines $p+1$ ``bands'' $\beta _i$
disjointly embedded in the interior of $F$ with the width of $\beta _i$ given by $w_i$, where $a_i$ traverses the length of $\beta
_i$, for $i=0,\ldots ,p$.  The bands coalesce near the boundary of $F$ and are attached one-to-the-next in the natural manner
to form ``one big band'' for each of the $\Delta$ distinguished points.  (More precisely, we are constructing here a ``partial
measured foliation from a weighted arc family regarded as a train track with 
stops''; for further details, see $\S$1 of [14]
or $\S$1.8 of [24].)

\vskip .1in

\noindent Letting $P$ denote the set of punctures of $F$, the {\it pure mapping class group} $PMC(F)$ is the group of isotopy
classes fixing $D\cup P$ pointwise  of all orientation-preserving 
diffeomorphisms of $F$ which fix $D\cup P$ pointwise.  

\vskip .1in

\noindent $PMC(F)$ acts on $Arc'(F)$ in the natural way,
and we define the {\it arc complex} of $F$ to be the 
quotient $$Arc(F)=Arc'(F)/PMC(F).$$  
Let us emphasize that distinct faces of a single simplex in $Arc'(F)$ may be
identified in this quotient, so $Arc(F)$ is a cell complex
but not a simplicial complex.  (See $\S$2.)  The complex described
in the Classical Fact is
$Arc'(F_{0,n}^0)=Arc(F_{0,n}^0)$ since $PMC(F_{0,n}^0)$ is trivial.

\vskip .1in

\noindent If $\alpha $ is an arc
family in $F$ with corresponding simplex $\sigma (\alpha )$ in $Arc'(F)$, then we shall let $[\alpha
]$ denote the
$PMC(F)$-orbit of $\alpha $ and $\sigma [\alpha ]$ denote the quotient of $\sigma (\alpha )$
in $Arc(F)$.  We may also abuse notation sometimes by suppressing the projective weight entirely and
write simply $[\alpha ]$ for the $PMC$-orbit of a projective positive weight on the arc family
$\alpha$.

\vskip .1in

\noindent Given two arc
families
$\alpha $ and
$\beta $ in $F$ whose union is again an arc family, the simplex $\sigma (\alpha \cup \beta )$ is
naturally identified with the simplicial join $\sigma (\alpha ) * \sigma (\beta )$.  In the same manner, if $F$
is the disjoint union of surfaces $F_1$ and $F_2$, then $Arc(F)$ is isomorphic to the simplicial join
$Arc(F_1)*Arc(F_2)$, and we use the 
convention that $K*\emptyset\approx K$, for any space $K$.

\vskip .1in

\noindent  
One difference between arc complexes in the bordered case from the
punctured case is that the
arcs in an arc family for a bordered surface come in a natural 
linear ordering which is invariant under the action
of the pure mapping class group.  

\vskip .1in

\noindent Specifically, enumerate once and for all the
distinguished points $p_1,p_2,\dots ,p_\Delta$ in the boundary of $F$ in any manner. 
The boundary of a regular neighborhood of 
$p_i$ in
$F$, for $i=1,2,\dots , \Delta$ contains an arc $A_i$ in the interior of $F$, and $A_i$ comes equipped with a
natural orientation (lying in the boundary of the component
of
$F-A_i$ which contains $p_i$).
If
$\alpha $ is an arc family in $F$, then there is a first intersection of a
component of $\alpha $ with the $A_i$ of least index, and this is the first arc $a_1$ in the
linear ordering.  By induction, the $(j+1)^{\rm st}$ arc (if such there be) in the putative linear
ordering is the first intersection (if any) of a component of $\alpha -\{ 
a_1,\dots ,a_{j-1}\}$ with the $A_i$ of least index. 

\vskip .1in

\noindent Since we take the pure mapping class group which fixes each distinguished point
and preserves the orientation of $F$, this linear ordering descends to a well-defined
linear ordering on the arcs in a
$PMC(F)$-orbit in $Arc'(F)$.  It follows immediately that the
quotient map $Arc'(F)\to Arc(F)$ is injective on the interior of any simplex in $Arc'(F)$, so that
the simplicial complex $Arc'(F)$ descends to a tractable
CW decomposition on $Arc(F)$.  (See $\S$2.) 

\vskip .1in 

\noindent Furthermore, in the
action of
$PMC(F)$ on
$Arc'(F)$, there can be no finite isotropy, and the isotropy subgroup of a simplex $\sigma (\alpha )$
in
$PMC(F)$ is either trivial or infinite.
The former case occurs if and only if each component of
$F-\cup \alpha $ is a polygon or a once-punctured polygon, 
and in this case we shall say that
$\alpha $ or $[\alpha ]$ {\it quasi fills} the surface $F$; for a compact surface $F$ with $s=0$, we shall simply
say that an arc family $\alpha$ or its orbit $[\alpha ]$ {\it fills} $F$ if each component of $F-\cup \alpha $ is a
polygon.

\vskip .1in

\noindent  In the extreme case that each component
of $F-\cup\alpha$ is a triangle or a once-punctured monogon, then
$\alpha$ is called a {\it quasi triangulation}, and in the compact case where components of the latter
type do not occur, $\alpha$ is called an {\it ideal triangulation} of $F$.  There are $N(F)+1$ arcs in a quasi triangulation of the bordered
surface $F$.

\vskip .2in

\noindent {\bf Lemma 2}~\it For $0\leq t<r$, the natural inclusion $F_{g,r}^s\subseteq F_{g,r-t}^{s+t}$ induces a
combinatorial equivalence
of $Arc(F_{g,r-t}^{s+t})$ with the subcomplex 
of $Arc(F_{g,r}^s)$ given by
$\{ \sigma [\alpha ]\subseteq
Arc(F_{g,r}^s):\alpha\cap\partial _i=\emptyset~{\rm for}~i=1,\ldots ,t\}$\rm

\vskip .2in

\noindent {\bf Proof}~~The natural inclusion $F_{g,r}^s\subseteq F_{g,r-t}^{s+t}$ induces an identification of
$Arc'(F_{g,r-t}^{s+t})$ with
$\{ \sigma (\alpha )\subseteq Arc'(F_{g,r}^s):\alpha\cap\partial _i=\emptyset~{\rm
for}~i=1,\ldots ,t\}$.  This identification is equivariant for the pure mapping class group actions
and induces
the asserted equivalence.~~~~\hfill{\it q.e.d.} 

\vskip .2in

\noindent{\bf Lemma 3}~\it There is a natural action of the $r$-torus $(S^1)^r$ on $Arc(F_{g,r}^s)$, and the corresponding
diagonal circle action is fixed point free.\rm

\vskip .2in

\noindent {\bf Proof}~The action of $(S^1)^r$ on $[\alpha ]\in Arc(F)$ for $F=F_{g,r}^s$  is defined in
terms of the band description of $[\alpha ]$ as follows.  For each $j=1,\ldots ,r$, consider the big band 
of $[\alpha ]$ incident on the $j^{\rm th}$
boundary component, and notice that at least one band must be non-empty.  Interpret $t_j\in S^1=[0,1]/(0\sim 1)$ and sweep a ratio
$t_j\in [0,1)$ of the
$j^{\rm th}$ big band around the
$j^{\rm th}$ boundary component, say to the right, to produce the resulting positive projectively weighted arc family under the action.

\vskip .1in

\noindent To see that the diagonal circle action is fixed point free, consider a positive weight $\mu _i$ on an arc family $\alpha=\{ a_i\}$ and take
the minimum
$$\mu _*={\rm min}\biggl \{{\rm min}~\{ \mu_i\} ,{\rm min}\bigr \{ \mu _i-\mu _j:\mu _i\neq \mu _j\bigr\}\biggr \}.$$
It is not difficult to see that $\mu ^{-1}_* ~\sum _i \mu _i$ is a lower bound for the order of the isotropy group in $S^1$ of $[\alpha ]\in
Arc (F)$.~~~~~\hfill {\it q.e.d.}

\vskip .2in

\noindent One of the key ingredients for proving that an arc complex
is non-spherical relies on the circle action in Lemma~3 and is illustrated
in the next result.  Letting $\approx$ denote PL-homeomorphism and $S^N$ the 
standard $N$-dimensional sphere, here is a specific example
of non-sphericity.

\vskip .2in

\noindent {\bf Proposition 4}~~\it It cannot be that both
$Arc(F_{0,2}^2)\approx S^5$ and $Arc(F_{0,1}^3)\approx S^3$. \rm 

\vskip .1in

\noindent (In fact, $Arc(F_{0,1}^3)\approx S^3$ as seen in Example~7 of 
Appendix~A, so $Arc(F_{0,2}^2)\not\approx S^5$.)  

\vskip .1in

\noindent{\bf Proof}~Suppose that both arc complexes are spherical.
According to Lemma~2, there are two copies of the 3-sphere $Arc(F_{0,1}^3)$ in the 5-sphere $Arc(F_{0,2}^2)$ corresponding
to those arc families that do not meet one or the other boundary component of $F_{0,2}^2$, and these 3-spheres are disjoint since each
non-empty arc must meet some boundary component.  Taking the quotient by the fixed point free diagonal circle action of Lemma~3, results
from algebraic topology [30] show that rationally the quotient of the 5-sphere is a complex projective
space
$CP^2$, and the homology class of each 3-sphere is non-trivial in the quotient.  Since the two 3-spheres are
disjoint in the manifold $Arc(F_{0,2}^2)$, their dual cohomology classes are supported on disjoint neighborhoods
and hence have trivial cup product.  This violates the ring structure of $CP^2$.~~~~\hfill{\it q.e.d.}

\vskip .2in

\noindent Given an arc family $\alpha $ in $F$, consider the surface $F-\cup\alpha$.  The non-smooth points or 
{\it cusps of} $\alpha$ on the frontiers of the components of $F-\cup\alpha $ give rise to distinguished points in the
boundary of
$F-\cup\alpha $ in the natural way, and this surface together with these distinguished points is
denoted
$F_{\alpha }$.  

\vskip .1in

\noindent We shall say that a surface $F$ is {\it inductive}  provided
$Arc(F_\alpha )$ is PL-homeomorphic to a
sphere of dimension
$N(F _\alpha)$ for each non-empty arc family
$\alpha $ in $F$.  

\vskip .1in
                                                                                
\noindent If $a$ is a separating arc whose endpoints coincide with the point $p$ in any surface,
then its two complementary components are distinguished by the fact that one component has exactly one
distinghished point
arising from $p$, and the other component has exactly two
distinguished points arising from $p$.  We
shall refer to these respective complementary regions as the {\it one-cusped}
 and
{\it two-cusped} components.

\vskip .3in

\noindent {\bf 2.  Barycentric subdivision of arc complexes}

\vskip .2in

\vskip .1in

\noindent This section is dedicated to proving the following result.

\vskip .2in

\noindent {\bf Lemma 5}~~\it If $F$ is inductive, then $Arc(F)$ is a PL-manifold of
dimension $N(F)$.\rm

\vskip .2in

\noindent The proof relies on the simplicial structure of the arc complexes, to which we now turn our attention.
Recall that a CW decomposition is said to be {\it regular} if the characteristic mapping of
each open cell extends to an embedding of the closed cell.

\vskip .2in

\noindent {\bf Lemma 6}~\it The first barycentric subdivision of $Arc'(F)$ descends to a regular CW
decomposition of $Arc(F)$.  Furthermore, the second barycentric subdivision of $Arc'(F)$ descends to
a simplicial complex on $Arc(F)$.\rm

\vskip .2in

\noindent {\bf Proof}~Let $\alpha _0\subseteq \alpha _1\subseteq \cdots \subseteq \alpha _p$ be a nested
collection of distinct arc families in $F$.  A $p$-simplex $\sigma$ in the first barycentric subdivision of
$Arc'(F)$ is given by a ``flag''
$$\sigma :~~ \sigma (\alpha _0)<\sigma (\alpha _1)<\cdots <\sigma (\alpha _p),$$
where $\sigma (\alpha )<\sigma (\beta )$ if $\sigma (\alpha )$ is a proper face of $\sigma (\beta )$, or in other
words,
$\alpha $ is a proper non-empty subset of the arc family $\beta $.  A $q$-dimensional face $\tau$ of $\sigma$ is
therefore given by a sub flag
$$\tau :~~\sigma (\alpha _{i_0})<\sigma (\alpha _{i_1})<\cdots <\sigma (\alpha _{i_q}),$$
where $i_0,i_1,\ldots ,i_q$ are distinct members of $\{ 0,1,\ldots ,p\}$.  Each $\alpha _i\subseteq \alpha
_p $ is evidently uniquely determined by its cardinality, so each face $\tau$ of $\sigma$ is likewise
uniquely determined by the tuple of cardinalities of its component arc families $\alpha _{i_0},\alpha
_{i_1},\ldots ,\alpha _{i_q}$.  Since such cardinality tuples are invariant under the action of $f\in PMC(F)$, if
$f(\tau )$ is a face $\tau '$ of
$\sigma$, then $\tau =\tau '$, so $f(\alpha _{i_q})=\alpha _{i_q}$.  Since $f$ respects the linear orderings on
component arcs in an arc family, $f$ must fix each arc in $\alpha _{i_q}$ proving the first assertion.

\vskip .1in

\noindent A $p$-simplex $\sigma$ in the second barycentric subdivision of $Arc'(F)$ is analogously given by a
``flag of flags''
$$\eqalign{
\sigma :~~&\sigma (\alpha _{0,0})<\sigma (\alpha _{0,1})<\cdots <\sigma(\alpha _{0,n_0})\cr
&\sigma (\alpha _{1,0})<\sigma (\alpha _{1,1})<\cdots <\sigma(\alpha _{1,n_1})\cr
&~~~~~\vdots ~~~~~~~~~~~~\vdots ~~~~~~~~~~~~~~~~~~~~\vdots\cr
&\sigma (\alpha _{p,0})<\sigma (\alpha _{p,1})<\cdots <\sigma(\alpha _{p,n_p})\cr
}$$
where $\{ \alpha _{i,0},\alpha _{i,1},\ldots ,\alpha _{i,n_i}\}$ is a proper non-empty subset of
$\{ \alpha _{i+1,0},\alpha _{i+1,1},\ldots ,\alpha _{i+1,n_{i+1}}\}$, for each $i=0,1,\ldots ,p-1$.  The maximal
entry
$\alpha _{i,n_i}$ in each row is called the ``carrier'' of the row, and the carrier
$\alpha _{i,n_i}$ is a not necessarily proper subset of
the consecutive carrier $\alpha _{i+1,n_{i+1}}$.  However, if the cardinality of a carrier $\alpha _{i,n_i}$ is
majorized by the cardinality of a carrier
$\alpha _{j,n_j}$, then $\alpha _{i,n_i}\subseteq \alpha _{j,n_j}$, for each $0\leq i,j\leq p$ (and though the
carrier
$\alpha _{i,n_i}$ is uniquely determined by its cardinality, the index $i$ of the carrier is not so
determined).

\vskip .1in
 \noindent Now, suppose that
$$\eqalign{
\sigma _1 :~~&\sigma (\alpha _{i_0,0})<\sigma (\alpha _{i_0,1})<\cdots <\sigma(\alpha _{i_0,n_{i_0}})\cr
&\sigma (\alpha _{i_1,0})<\sigma (\alpha _{i_1,1})<\cdots <\sigma(\alpha _{i_1,n_{i_1}})\cr
&~~~~~~\vdots ~~~~~~~~~~~~~\vdots ~~~~~~~~~~~~~~~~~~~~~~~~\vdots\cr
&\sigma (\alpha _{i_\mu ,0})<\sigma (\alpha _{i_\mu ,1})<\cdots <\sigma(\alpha _{i_\mu ,n_{i_\mu }})\cr
}~~~~\eqalign{
\sigma _2 :~~&\sigma (\alpha _{j_0,0})<\sigma (\alpha _{j_0,1})<\cdots <\sigma(\alpha _{j_0,n_{j_0}})\cr
&\sigma (\alpha _{j_1,0})<\sigma (\alpha _{j_1,1})<\cdots <\sigma(\alpha _{j_1,n_{j_1}})\cr
&~~~~~~\vdots ~~~~~~~~~~~~~\vdots ~~~~~~~~~~~~~~~~~~~~~\vdots\cr
&\sigma (\alpha _{j_\nu ,0})<\sigma (\alpha _{j_\nu ,1})<\cdots <\sigma(\alpha _{j_\nu ,n_{j_\nu }})\cr
}$$
are two faces of $\sigma$, where $\nu\geq\mu$, and $f_k\in PMC(F)$ maps $\sigma _k$, for $k=1,2$, to faces of a
common simplex $\tau$.  Since $\alpha _{i_{\mu},n_{i_\mu}}\subseteq \alpha _{j_\nu,n_{j_\nu}}$, the cardinality
of 
$f_1(\alpha _{i_{\mu},n_{i_\mu}})$ is majorized by the cardinality of $f_2(\alpha _{j_\nu,n_{j_\nu}})$, so 
$f_1(\alpha _{i_\mu,n_{i_\mu}})\subseteq f_2(\alpha _{j_\nu,n_{j_\nu}})$ by the earlier remarks.  Thus,
$f_2^{-1}\circ f_1(\alpha _{i_\mu , n_{i_\mu}})=\alpha _{i_\mu ,n_{i_\mu}}$, so $f_2^{-1}\circ f_1$ fixes each
arc of
$\alpha _{i_\mu ,n_{i_\mu}}$ as before, and $\sigma$ and $\tau$ in fact meet along the face of $\sigma$
corresponding to the union of the rows $i_0,i_1,\ldots ,i_\mu$ and $j_0,j_1,\ldots j_\nu$, proving the second
assertion.~~~~~\hfill{\it q.e.d.}

\vskip .2in

\noindent Suppose that $\alpha =\{ a_0,a_1,\ldots ,a_p\}$ is an arc family in $F$ and $\beta =\{ b_0,b_1,\ldots
,b_q\}$ is an arc family in $F_\alpha$, where the component arcs are in each case
enumerated in their canonical linear orderings.  For any permutations
$ i_0,i_1,\ldots ,i_p$ of $\{ 0,1,\ldots ,p\}$ and $j_0,j_1,\ldots ,j_q$
of $\{ 0,1,\ldots ,q\}$, there is an associated $p$-dimensional simplex
$$
\sigma :~~\sigma (\{a_{i_0}\}) <\sigma (\{ a_{i_0},a_{i_1}\})<\cdots <\sigma (\alpha )$$
in the first barycentric division of $Arc'(F)$, an associated $q$-dimensional
simplex $$\tau :~~\sigma (\{ b_{j_0}\} )
<\sigma (\{ b_{j_0},b_{j_1}\} )<\cdots <\sigma (\beta)$$
in the first barycentric subdivision of $Arc'(F_\alpha )$, and a 
corresponding
$(p+q)$-dimensional simplex
$$\rho:~\sigma (\{a_{i_0}\}) <\sigma (\{ a_{i_0},a_{i_1}\})<\cdots <\sigma (\alpha )<\sigma (\alpha \cup\{
b_{j_0}\} ) <\sigma (\alpha \cup\{ b_{j_0},b_{j_1}\} )<\cdots <\sigma (\alpha\cup\beta)$$
in the first barycentric subdivision of $Arc'(F)$.
Given the respective affine coordinates $[u_0:u_1:\cdots :u_p]$ and $[v_0:v_1:\cdots :v_q]$ for projective
weightings on $\sigma$ and $\tau$ and given any $0\leq t\leq 1$, there is an induced  projective weighting on
$\rho$ with affine coordinates $$[(1-t)u_0:(1-t)u_1:\cdots :(1-t)u_p:tv_0:tv_1:\cdots :tv_q].$$
Varying over all pairs of permutations, this correspondence determines a continuous mapping $\sigma (\alpha) *
Arc'(F_\alpha )\to Arc'(F)$, which descends to a well-defined continuous mapping
$$\sigma [\alpha ] * Arc(F_\alpha )\to Arc(F).$$

\vskip .2in

\noindent {\bf Lemma 7}~\it For any arc family $\alpha $ in $F$, the natural mapping 
$\sigma [\alpha ] * Arc(F_\alpha )\to Arc(F)$
is an embedding on $Int~\sigma [\alpha ] *Arc(F_\alpha )$.\rm

\vskip .2in

\noindent{\bf Proof of Lemma~7}~Suppose that $\alpha ^\ell =\{a_0^\ell,a_1^\ell ,\ldots ,a_p^\ell\}$ is an arc
family in $F$, $\beta ^\ell =\{ b_0^\ell ,b_1^\ell ,\ldots ,b_q^\ell \}$ is an
arc family in
$F_{\alpha ^\ell}$, $ i_0^\ell ,i_1^\ell ,\ldots ,i_p^\ell $ is a permutation of $\{ 0,1,\ldots ,p\}$, and
$j_0^\ell ,j_1^\ell ,\ldots ,j_q^\ell $ is a permutation of $\{ 0,1,\ldots ,q\}$ with corresponding simplices
$\rho ^\ell$ in the first barycentric subdivision of $Arc'(F)$, for $\ell =1,2$, 
where $f(\rho ^1)=\rho ^2$ for some $f\in PMC(F)$.  Since $f$ preserves cardinalities, we must have $f(\alpha
^1)=\alpha ^2$.   Since $f$ respects the linear orderings on the arcs comprising an arc family, we must
furthermore have $f(a_m^1)=a_m^2$, whence $i_m^1=i_m^2$, for each $m=0,1,\ldots ,p$.  In precisely the same
manner, we conclude that $f(b_n^1)=b_n^2$ and $j_n^1=j_n^2$, for each $n=1,2,\ldots ,q$, and the result
follows.~~~~\hfill{\it q.e.d.}

\vskip .2in

\noindent{\bf Proof of Lemma 5}~We shall prove that any point $x\in Arc(F)$ admits a neighborhood
which is PL-homeomorphic to an open ball of dimension $N=N(F)$.  Of course, $x$ lies in the interior $Int~\sigma
[\alpha ]$ of $\sigma [\alpha ]$, for some arc family $\alpha$ in $F$, and we suppose that $\alpha$ is comprised
of
$p+1$ component arcs. We have the identity $N+1= (N(F_\alpha )+1)+(p+1)$ since both sides of the equation give the
number of arcs in a quasi triangulation of $F$.  Since $Arc(F_\alpha )\approx S^{N(F_\alpha
)}$ by hypothesis, $Int~ \sigma [\alpha ]*Arc(F_\alpha )$ is PL-homeomorphic to an open ball of
dimension
$p+N(F_\alpha )+1 =N$, which by the previous lemma gives the required neighborhood of $x$ in $Arc(F)$.
~~~~\hfill{\it q.e.d.}

\vskip .3in

\noindent {\bf 3. The cusp operator and Classical Fact}

\vskip .2in

\noindent Given a surface $F$, suppose that $F'$ arises from $F$ by adding one cusp on the boundary, 
i.e., choosing
one extra distinguished point on the boundary of $F$.  We shall think of  a cusp $p\in\partial F$
as splitting into two cusps $q_0,q_1$ of $F'$ as depicted in Figure~1a, where an embedding of $F$
in $F'$ is depicted; 
the natural identification of the other cusps
of $F$ with those of $F'$ is also indicated.  We may identify $S^0=\{ a_0,a_1\}$, where the arcs
$a_0,a_1$ are as illustrated in Figure 1b.

\hskip .3in {{{\epsffile{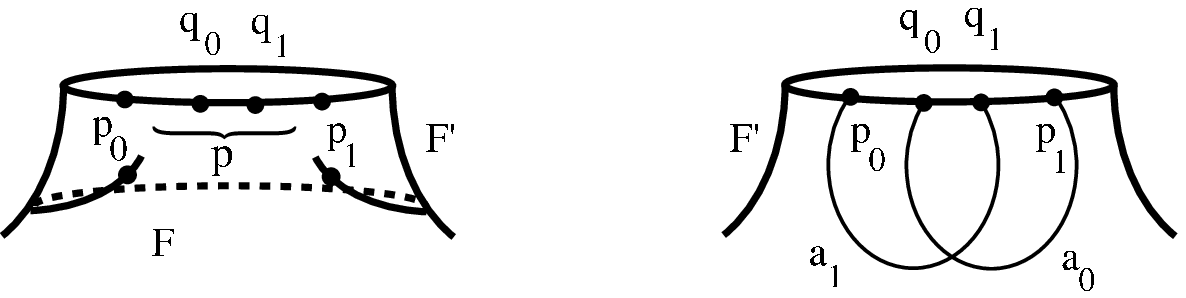}}}}

\hskip .2in\noindent {\bf Figure 1a}~The embedding of $F$ in $F'$.\hskip .5in
{\bf Figure 1b}~The arcs $a_0,a_1$ in $F'$.

\vskip .1in

\centerline{{\bf Figure 1}~The embedding $F\subseteq F'$ and special arcs in $F'$.}

\vskip .1in

\noindent A point of
$S^0*Arc(F)$ is thus given by the projective class of~a weight $0\leq w\leq 1$ on $a_0$ or $a_1$ and a non-negative $\nu$ on
some arc family $\beta$ in $M$, where at least one of $w,\nu$ is non-zero.  Let $n\geq 0$ denote the total weight of all arcs in $\beta$
which meet
$p$, i.e., the width of the one big band of $(\beta ,\nu )$ at $p\in F$.

\vskip .1in

\noindent Now we take the blow-up of $S^0*Arc(F)$ along the codimension-one sphere $Arc(F)$ by an interval in the join direction, where
this interval is identified with the one-simplex of all
probability measures $\rho$ on the set $\{ 0,1\}$.  This defines a blow-up
$\bigl (S^0*Arc(F)\bigr )^\wedge$, where a point of this blow-up is given by the projective class of $\nu ,w$ together
with a probability density on
$\{ 0,1\}$ whenever $w=0$.

\vskip .1in

\noindent Define the ``cusp operator''
$$\hat\Phi : \bigl (S^0*Arc(F)\bigr )^\wedge ~~\to ~~
Arc(F')$$
as follows.  If $w\neq 0$ is a weight on the arc $a_i$, for $i=0,1$, then
connect the arcs of $\beta$ 
hitting
$p$, if any, to $q_i$, as illustrated in Figure 2a.  If $w=0$, then use the probability density
$\rho$ on
$\{ 0,1\}$, to split the big band and send proportion $\rho (i)$ of the arcs of $\beta$ hitting $p$, if any, to
$q_i$, for
$i=0,1$, as illustrated in Figure 2b.  

\vskip .1in
 \noindent This defines a projectively weighted arc family in $F'$ as required to complete the definition of
$\hat\Phi$.
Notice that if
$\beta$ does not hit
$p$, then
$\rho$ plays no role in the construction, so the entire interval of probability
densities is collapsed to a single point. 

\vskip .05in

\centerline{{{\epsffile{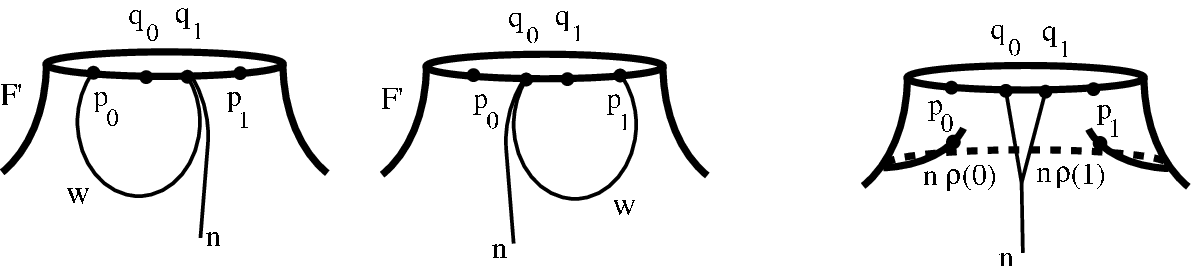}}}}

\hskip .2in{\bf Figure 2a} Extension for $w\neq 0$.\hskip 1.in{\bf Figure
2b} Extension for $w=0$.

\vskip .1in

\centerline{{\bf Figure 2}~Extending curves from $F$ to $F'$ for the cusp operator.}

\vskip .2in

\noindent {\bf Lemma 8}~ \it Suppose that $F'$ arises from $F$ by adding one cusp and that $F$ is inductive.  Then $Arc(F')$
is PL-homeomorphic to the suspension of $Arc(F)$.~~~~~\hfill \rm

\vskip .2in

\noindent {\bf Proof} 
The mapping $\hat\Phi$ is continuous by construction and is easily shown to be surjective.
$\hat\Phi$ is also easily shown to be injective
over the open set $$\{ Int~\sigma [\alpha ]\subseteq Arc(F'):\alpha~{\rm meets~every~distinguished~point~of}~F'\} .$$
We may consider the
quotient $\bigl (S^0*Arc(F)\bigr )^\wedge /\sim$, where $\sim$ collapses to a point each probability interval whenever the underlying arc
family in
$F$ does not meet $p$.  It follows that $\hat\Phi$ induces a continuous bijection $$\bar\Phi ~:~~\bigl (S^0*Arc(F)\bigr )^\wedge
/\sim ~~\to ~~Arc(F').$$ 

\vskip .2in

\noindent We may explicitly construct the continuous inverse to $\bar\Phi$ as follows.  First note that the subspace of arcs not meeting
a given cusp or not containing a given arc class is a full sub-complex of the arc complex.
In the notation of Figure~1, consider the closed subsets
$A_i=\{ \sigma [\alpha ]\subseteq Arc(F'): \alpha~{\rm does~not~meet}~q_i\}$, for $i=0,1$, 
and $B=\{\sigma [\alpha ]\subseteq Arc(F'):\{a_0,a_1\}\cap\alpha=\emptyset\}$ which form a cover
of $Arc(F')$.  Identify the subspace of $Arc(F')$ comprised of arc families that do not meet $q_0$ or 
do not meet $q_1$ with the subspace 
in the domain of $\bar\Phi$ corresponding to probability density $\rho$ concentrated at 0 or 1, each component of which is
PL homeomorphic to $Arc(F)$.  Invert the construction
in Figure 2a in the natural way to further define the inverse on  $A_0\cup A_1$, and likewise invert the construction in Figure 2b in the natural way to define
the inverse on
$B$.  

\vskip .2in

\noindent
Thus, $\bar\Phi$ is a PL-homeomorphism.  Results from algebraic topology [31]
show that $\bigl (S^0*Arc(F)\bigr )^\wedge /\sim$ and $\bigl (S^0*Arc(F)\bigr )^\wedge$
are PL homeomorphic provided $Arc(F)$ is a manifold, as follows from Lemma~5 and the hypothesis that $F$ is inductive.  
Furthermore, $\bar\Phi$ is a PL
homeomorphism
$Arc(F')\approx \bigl (S^0*Arc(F)\bigr )^\wedge /\sim$ as above, and blowing up and down finally gives
$\bigl (S^0*Arc(F)\bigr
)^\wedge\approx S^0*Arc(F)$, completing the proof.~~~~\hfill{\it q.e.d.}

\vskip .2in
  
\noindent {\bf Proof of the Classical Fact}~The proof 
follows by induction on the number of the sides of the
polygon using Lemma~8 with basis step given by the quadrilateral, which is treated in Example~1 of 
Appendix~A.~~~~~\hfill{\it q.e.d.}

\vskip .1in

\noindent A different proof of the Classical Fact (which
uses Thurston theory) is given in [25].

\vskip .3in

\noindent {\bf 4. The puncture operator, multiply punctured polygons, and $F_{1,1}^1$}

\vskip .2in

\noindent Suppose that $F=F_{g,r}^s$ contains a simple closed curve that separates a surface $F_{0,2}^1$ from a surface
$M=F_{g,r}^{s-1}$, where $F_{0,2}^1$ contains a fixed puncture of $F$.  Let $T=F_{0,1}^2$, so
$Arc(T)\approx S^1$ is analyzed in Example~4 of Appendix A, whose notation we adopt here.  We shall first define a real blow-up of 
$Arc(T)*Arc(M)$ under the assumption that $Arc(M)$ and $Arc(T)$ are PL-homeomorphic to spheres.  Thus, $Arc(T)*Arc(M)$ is again PL-homeomorphic
to a sphere, which we blow up along a codimension-one sphere, in order to define the puncture operator. 

\vskip .1in

\noindent Suppose that $m\in Arc(M)$, and consider the
2-disk
$$D^2=Arc(T)*\{ m\}\approx S^1*\{ m\}\subseteq Arc(T)*Arc(M).$$  As illustrated in Figure~3, the two singleton arc classes
$[a_0],[a_1]$ in
$T$ are taken to lie on the real axis, $m\in Arc(M)$ is taken as the origin, and $Arc(T)$ is taken as the unit circle in ${\bf R}^2$.  There
is a one-dimensional blow-up of the real unit interval in
$D^2$ illustrated in Figure~3, and we let
$\bigl (Arc(T)*Arc(M)\bigr )^\wedge$ denote the corresponding one-dimensional blow-up 
of $Arc(T)*Arc(M)\approx S^{N(F)}$ along the
codimension-one sphere $\{ [a_0],[a_1]\}*Arc(M)$.  

\vskip .1in

\centerline{\epsffile{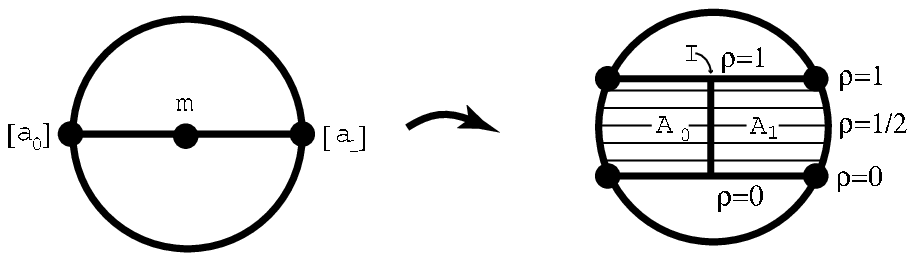}}

\vskip .1in

\centerline{{\bf Figure 3}~The fundamental blow-up.}

\vskip .2in

\noindent Figure~3 furthermore indicates the level sets of the blow-up parameter
$0\leq
\rho\leq 1$, where $\rho\equiv 0$ extends below and $\rho\equiv 1$ above the real axis to a continuous function
$D^2\to [0,1]$, and indeed to a continuous function
$$\rho:\bigl (Arc(T)*Arc(M)\bigr )^\wedge \to [0,1].$$
In particular, $m\in Arc(M)$ is blown-up to an interval $I$ in $D^2$, which decomposes the blown-up locus in $D^2$ into rectangles
$A_0,A_1$, as is also illustrated in this figure.

\vskip .2in

\noindent 
A point $\zeta\in (Arc(T)*Arc(M))^\wedge$ is given by the simultaneous projectivization
of a pair of weighted arc families $(\alpha,\mu)$ in $T$ and $(\beta ,\nu)$ in $M$, where at least one of $\mu$
or $\nu$ is non-vanishing; on $\{ [a_0],[a_1]\} *Arc(M)$, there is a further continuous blow-up
parameter
$0\leq
\rho\leq 1$.

\vskip .1in

\noindent Define the ``puncture operator'' $$\hat\Phi :(Arc(T)*Arc(M))^\wedge \to Arc(F)$$ as follows.  
Identify $Arc(T)$ with the subset of
$Arc(F_{0,2}^1)$ corresponding to arc families that do not meet a fixed boundary component $\partial _*$ of $F_{0,2}^1$ using Lemma~2, and let
$\partial ^*=\partial F_{0,2}^1-\partial _*$.   Identify a boundary component of 
$M$ with
$\partial _*$ to get a surface homeomorphic to $F$ and let $n$ denote the width of the one big band of $(\beta ,\nu )$
at $\partial _*$.  

\vskip .1in

\noindent If $n=0$, then $\hat\Phi (\zeta )$ is defined to be the projective class of the disjoint union of $(\alpha ,\mu )$ and $(\beta
,\nu )$ in
$F$.  

\vskip .1in

\noindent In the generic case that $\partial _*$ lies in a monogon complementary to $\alpha$,  
there
is a unique way to extend
$(\beta ,\nu )$ to a weighted arc family in the complement of
$\alpha$ by running a band of width $n$ across the monogon complementary to $\alpha$; combining $(\alpha ,\mu )$ and $(\beta ,\nu )$ in
this manner gives a weighted arc family $(\gamma ,\lambda )$ in $F$.  Taking the projective class $\hat\Phi (\zeta)$ of $(\gamma
,\lambda)$ defines the puncture operator on the complement of the blow-up, i.e., on the complement
of the  union of the rectangles $A_0,A_1$ in Figure~3.

\vskip .1in 
\noindent On the rectangle $A_1$, there is furthermore a one-parameter family
$(\gamma _\rho,\lambda _\rho )$, for $0\leq \rho\leq 1$, as follows: sweep to the right around $\partial ^*$ a proportion $\rho$ of the big band
at $\partial ^*$ of $(\gamma _0,\lambda _0)$ defined above.  This defines $\hat\Phi$ on $A_1$.

\vskip .1in

\noindent In the remaining non-generic case that $\partial _*$ lies in a bigon complementary to $\alpha$ (i.e. on $A_0-I$), we may split the one
big band at $\partial _*$ into two bands that run to the two distinguished points in the two possible ways in this bigon; there is thus again a one-parameter
family
$(\gamma _\rho,\lambda _\rho )$, for $0\leq \rho\leq 1$, where proportion $\rho$ of $(\gamma _0,\lambda _0 )$ goes to the right and $1-\rho$ to the left
along $\partial ^*$ in this bigon.

\vskip .2in

\noindent {\bf Lemma 9}~\it ~The puncture operator $\hat\Phi$ is continuous.  
Provided $F$ has a single boundary component, $F$ is inductive, and
$Arc(M)$ a sphere, the puncture operator $\hat\Phi$ furthermore has degree one.
\rm

\vskip .2in

\noindent {\bf Proof}~~Continuity of $\hat\Phi$ follows from consideration of
Example~4 of Appendix~A, and in particular from Figure~A.1, which illustrates $Arc(T)$.
In effect, the construction defining $\hat\Phi$ 
on $A_0$ corresponds to sweeping proportion
$\rho$ of the big band at $\partial _*$ for $(\beta ,\nu )$ from the left to the right around $\partial ^*$, while on $A_1$, it corresponds to twisting to
the right proportion
$\rho$ of the the big band at $\partial ^*$.  These two constructions agree on the common interface $I$, i.e., the
sweeping agrees with the twisting if there is no arc in $T$.  
It is furthermore clear from this discussion that the
construction likewise extends continuously across the top and bottom of 
$A_0,A_1$, and therefore $\hat\Phi$ is indeed continuous.  

\vskip .1in

\noindent Since $Arc(M)$ is a sphere by hypothesis,  so too $Arc(T)*Arc(M)\approx S^1* Arc(M)$ is again a sphere.  If $Arc(F)$ is inductive, then $Arc(F)$ is a
manifold (by Lemma~5), so it makes sense to discuss the degree of $\hat\Phi$.  We finally show that $\hat\Phi$ has degree one provided $F$ (and hence $M$) has
only one boundary component.

\vskip .1in

\noindent To this end, suppose that $(\gamma ,\lambda)$ is a weighted arc family in $F$, and let $b$ be an arc in $F$
bounding a monogon which contains a puncture $p$ of $F$.  Perhaps $p$ lies in a monogon complementary to $\gamma$ in which case $\gamma$ contains
an element of the $PMC(F)$-orbit $[b]$ of $b$. 
If not, then $p$ lies in complementary component to $\gamma$ which has smaller Euler characteristic than the once-punctured monogon, and there is an arc
bounding a monogon containing $p$ in this component.  Since there is only one boundary component whose unique cusp gives rise to all cusps of
$F_\gamma$, this arc lies in the class $[b]$.
Thus, if $\gamma$ does not contain an
element of
$[b]$, then we may add an element of
$[b]$ complementary to
$\gamma$.  Furthermore, if $\gamma$ contains an element $b$ of $[b]$, then it contains a unique such element
since $b$ is separating.  There is thus a continuous surjection $\{ [b]\} * Arc(F_{\{ b\} })\to Arc(F)$
which is an embedding on the interior
when $F$ has a unique boundary component with a single distinguished point.
Notice that $\{ [b]\} * Arc(F_{\{ b\} })$ is a PL ball of dimension $N(F)$ since $F$ is inductive.

\vskip .1in

\noindent Taking $p$ to be the added puncture of the puncture operator and $b=a_0$, it follows that there is a continuous surjection
$\{ [a_0]\}* Arc(M') \to Arc(F)$,
where $M'$ of type $M_{g,(2)}^{s-1}$ has two distinguished points on the boundary
whereas $M$ of type $M_{g,1}^s$ has only one.  Moreover, this mapping is an embedding on the interior.

\vskip .1in

\noindent This defines the horizontal mapping on the bottom of the diagram

\vskip .2in

\settabs 7\columns

\+&$(Arc(T)*Arc(M))^\wedge$&&$\leftrightarrow$&$Arc(T)*Arc(M)$\cr
\+\cr
\+&~~~~~~~$\hat\Phi~~\downarrow$&&&~~~~~~~~~~$\uparrow$\cr
\+\cr
\+&~~~~~~~$Arc(F)$&&$\leftarrow$&$\{ [a_0]\} *Arc(M'),$\cr

\vskip .2in

\noindent and the horizontal mapping on the top is the homotopy equivalence given by blowing up and down. 
Finally making the identifications $Arc(M')\approx S^0*Arc(M)\approx \{b_0,b_1\}*Arc(M)$ according to Lemma~8, 
we may define the surjective mapping on the right sending $\{ b_i\} *Arc(M)$ to $\{ a_1\}*Arc(M)$, for $i=0,1$,
which is an embedding on the interior.  Since this diagram of continuous maps commutes (starting from the lower-right corner) on a dense set,
it commutes up to homotopy.  Since each map other than
$\hat\Phi$ has degree one, $\hat\Phi$ itself 
thus also has degree one. 
~~~~~\hfill{\it q.e.d.}

\vskip .2in

\noindent {\bf Corollary 10}~~\it In the notation of the puncture operator, 
assume that $F$ is inductive with $Arc(F)$ a manifold of dimension at least five
and that $M$ has exactly one boundary component. 
If $Arc(M)$ is spherical, then $Arc(F)$ is spherical. 

\rm

\vskip .2in

\noindent {\bf Proof}~According to Lemma~9 and assuming that
$Arc(F)$ is a manifold, we have a continuous degree one PL map $\hat\Phi :S^N\to Arc(F)$
from the
$N=N(F)$-sphere
$S^N$ to the $N$-dimensional manifold $Arc(F)$.
Smale's h-cobordism theorem [29] shows that provided $N\geq 5$, $Arc(F)$ is PL-homeomorphic to
$S^N$.~~~~\hfill{\it q.e.d.}

\vskip .2in

\noindent {\bf Proposition 11}~\it $Arc(F_{1,1}^1)$ is spherical.\rm

\vskip .2in

\noindent {\bf Proof}~We have $Arc(F_{0,1}^2)*Arc(F_{1,1}^0)\approx S^1*S^3$ using Examples~4 and 6 of Appendix~A.  Since $F=F_{1,1}^1$ is
inductive by Proposition~A, it follows that $Arc(F)$ is a manifold by Lemma~5.  It follows from Corollary~10 that
$Arc(F)$ is spherical. ~~~~\hfill{ q.e.d.}

\vskip .2in

\noindent {\bf Proposition 12}~\it The arc complex of any multiply punctured polygon is spherical.\rm

\vskip .2in

\noindent {\bf Proof}~~Use Lemma~8, the Classical Fact, and strong induction on
$s\geq 4$ to conclude that $F_{0,1}^{s+1}$ is inductive.   By Lemma~5, $Arc(F_{0,1}^{s+1})$ is a
manifold and hence is spherical.~~~~\hfill{ q.e.d.}

\vskip .3in

\noindent {\bf 5. Counter-examples}

\vskip .2in

\noindent Together with the appendices, we have proved sphericity of all the arc complexes asserted to be spherical in Theorem~1, and it
remains to prove non sphericity of the other arc complexes.  To this end, we have the following result: 

\vskip .2in

\noindent {\bf Lemma 13}\it ~~Suppose that $\alpha$ is an arc family in $F$ and $M$ is a component of $F_\alpha=F-\cup\alpha$.  If $Arc(M)$ is
not spherical, then $Arc(F)$ not a PL-manifold.\rm

\vskip .2in

\noindent {\bf Proof}~~By Lemma~7, the link of the simplex $\sigma [\alpha ]$ in the second barycentric subdivision is not
spherical, so no point of $Int~\sigma[\alpha ]$ is a manifold point of $Arc(F)$. ~~~~\hfill{\it q.e.d.}

\vskip .2in

\noindent We shall say that the surfaces $F_{0,4}^0$, $F_{0,3}^1$, $F_{0,2}^2$, $F_{1,2}^0$ are {\it type 1 surfaces}.
As we shall see, they are the only inductive surfaces with non spherical arc complexes.

\vskip .2in

\noindent {\bf Lemma 14}~~\it The type 1 surfaces have non spherical arc
complexes.\rm

\vskip .2in

\noindent {\bf Proof}~~Enumerate
the boundary components
$\partial _1,\ldots ,\partial _r$ of
$F$, and let
$S_I\subseteq Arc(F)$ denote the geometric realization of the set of all arc families in $F$ which are disjoint from $\partial _i$ for each
$i\in I$.  Thus,
$S_I\cap S_J=S_{I\cup J}$ and
$S_{\{ 1,\ldots ,r\}}=\emptyset$.

\vskip .1in

\noindent For each $F_{0,r}^s$ with $r+s=4$ and each $i=1,\ldots ,r$, we have $S_{\{ 1,\ldots ,r\} -\{ i\}}\approx Arc(F_{0,1}^3)\approx
S^3$ by Lemma~2 and Example~7 of Appendix~A.

\vskip .1in

\noindent Since $S_{\{ 1\}
}\cap S_{\{ 2\} }=\emptyset$ for $F_{0,2}^2$, the argument in Proposition~4 shows
that the arc complex of $F_{0,2}^2$ is not spherical.  Likewise, since $S_{\{ 1,2\}
}\cap S_{\{ 1,3\} }=\emptyset$ for $F_{0,3}^1$ and since $S_{\{ 1,2,3\} }\cap S_{\{ 1,2,4\} }=\emptyset$
for $F_{0,4}^0$, $Arc(F_{0,3}^1)$ and $Arc(F_{0,4}^0)$ are also non spherical.

\vskip .1in

\noindent In the remaining case $F_{1,2}^0$, define $S_{\{ i\} }\subseteq Arc(F_{1,2}^0)$
as before, for $i=1,2$.  Again, Lemma~2 and Proposition~11 imply that $S_{\{ 1\} }\approx S_{\{ 2\} }\approx Arc (F_{1,1}^1)\approx S^5$, 
where $S_{\{ 1\}
}\cap S_{\{ 2\} }=\emptyset$, so again the
argument of Proposition~4 shows that $Arc(F_{1,2}^0)$ is non spherical.~~~~~\hfill{\it q.e.d.}

\vskip .2in

\noindent {\bf Lemma 15}~\it ~For any surface $F$ other than those 
already proved to be spherical, we have $M\subseteq F$ for one or
more of the type 1 surfaces $M$.\rm

\vskip .2in

\noindent {\bf Proof:}~If $g\geq 3$ or $g=2, r>1$ or $g=1, r>2$,
then $F_{0,4}^0\subseteq F$; if $g=2, s>1$, then $F_{0,3}^1\subseteq F$; if $g=1$ and $s>1$, then $F_{0,2}^2\subseteq F$; and if $g=2$ and
$r=1$, then $F_{1,2}^0\subseteq F$.~~~~~\hfill{\it q.e.d.}

\vskip .2in

\noindent {\bf Corollary 16}~\it The type 1 surfaces are the unique inductive
surfaces with non spherical arc complexes.\rm

\vskip .2in

\noindent {\bf Proof}~~In light of the two previous lemmas, it remains 
only to prove that the type 1 surfaces are
inductive.  This follows from the positive sphericity results proved before.~~~~~\hfill{\it q.e.d.}

\vskip .3in

\noindent {\bf 6. Proof of Theorem~1 and concluding remarks}

\vskip .2in

\noindent {\bf Proof of Theorem~1}~We have proved that all
the asserted surfaces do indeed have spherical arc complexes in the appendices and in Propositions~11 and 12 plus the Classical Fact using
Lemma~8.  By Lemmas~13-16, no other surface has a spherical arc complex, and 
furthermore, only the type 1 surfaces have manifold
arc complexes.~~~~~\hfill{\it q.e.d.}

\vskip .2in

\noindent It is worth explicitly pointing out here that the examples in Appendix~A show that the link of any simplex in the second
barycentric subdivision of $Arc(F)$ of codimension at most four is a sphere, and Theorem~1 also gives a complete classification of the manifold
points of
$Arc(F)$ for any bordered surface $F$.  

\vskip .1in

\noindent Arc complexes as piecewise linear (PL) objects conjecturally have specific
singularities and topology which are described recursively as follows.  
Say that a PL sphere is a type zero space.  The four specific arc complexes of type 1 surfaces as PL objects are the type one
spaces.  Inductively for $n>1$, define a
type
$n$ space to be a finite polyhedron, defined up to PL-isomorphism, so that the link of each vertex in any compatible
triangulation is PL isomorphic to an iterated suspension of  the join of at most two spaces of type  less than
$n$.  We conjecture that each arc complex of a connected bordered surface determines an underlying space of finite
type $n<\infty$.  Indeed, if one could drop the hypothesis that $F$ is inductive in Lemma~8, then the conjecture would
follow as in the case of Theorem~1.  Perhaps a recursive application of the techniques of [31] in the context of the inductive proof
of Theorem~1 could provide the required collapse and describe the general stratified structure of arc complexes as
governed by the four type 1 surfaces.

\vskip .1in

\noindent Suppose that $F$ is a possibly bordered surface, and let ${\cal PF}_0={\cal PF}_0(F)$ denote Thurston's space of projective
foliations 
of compact support in $F$.  It is well-known (cf. the Remark after Lemma~B or 
[5,24,27]) that ${\cal PF}_0$ is a PL sphere of dimension
$6g-7+2(r+s)$.  According to Thurston theory [5,27], there is a compactification $\overline T=\overline T(F)$ of the Teichm\"uller space
$T=T(F)$ of
$F$, where ${\cal PF}_0=\overline T-T$, $\overline T$ is a closed ball of dimension $6g-6+2(r+s)$, and the action of $PMC=PMC(F)$ on
$T$ extends continuously to the natural action on $\overline T$.  
As is also well-known (cf. [24]), simple closed curves are dense in
$P{\cal F}_0$, so the quotient
${\cal PF}_0/PMC$ is maximally non Hausdorff in the sense that its largest Hausdorff quotient is a singleton.  

\vskip .1in

\noindent Likewise, we may
consider the space ${\cal PF}(F)$ of all projective measured foliations, not necessarily of compact support on $F$, which is again 
seen to be a sphere of dimension $6g-7+3r+2s+\Delta$ (cf. the 
Remark after Lemma B or [5,24]); for the same reasons as before, the
quotient ${\cal PF}/PMC\supseteq {\cal PF}_0/PMC$ is maximally non Hausdorff.  To our knowledge, not much has been done to study these
quotients, and the point we wish to make here is that the non Hausdorff space ${\cal PF}(F)/PMC(F)$ contains $Arc(F)$ as an open dense set,
and $Arc(F)$ is not only Hausdorff but also a stratified space of a particular type.

\vskip .1in

\noindent  As in [21,25] and as we shall discuss in a forthcoming paper, the strictly combinatorial statement of Theorem~1 has immediate
interpretations in terms of the combinatorics of macromolecular folding in computational biology.  Furthermore, the operadic structures
on arc complexes and Riemann moduli spaces [2,3,4,14,28,32] might be studied using the combinatorial structure described here, and indeed conversely as well.

\vskip .1in

\noindent Turning finally to the relation between the sphericity 
results and the homology of Riemann's moduli space, let us for
simplicity restrict
attention to the case $s=0$ of compact surfaces.  (More generally, one works with the one-point compactification of a punctured surface.)
Let $\partial $ denote the boundary mapping of the chain
complex
$\{ C_p (Arc):p\geq 0\}
$ of  $Arc=Arc(F)$.   Suppose that
$\alpha$ is an arc family in
$F$ with corresponding cell $\sigma [\alpha ]\in C_p(Arc)$.  A codimension-one face of
$\sigma [\alpha ]$ of course corresponds to removing one arc from $\alpha$, and
there is a dichotomy on such faces $\sigma [\beta ]$ depending upon whether the rank
of the first homology of 
$F_\beta$ agrees with that of $F_\alpha$ or differs by one from that of $F_\alpha$.
This dichotomy decomposes $\partial $ into the sum of two operators $$\partial =\partial
^1+\partial ^2,$$ where $\partial ^2$ corresponds to the latter case.

\vskip .1in

\noindent The operators
$\partial ^1,\partial ^2$ are a pair of anti-commuting differentials,
so there is a  spectral sequence
converging to $H_*(Arc)$ corresponding to the bi-grading
$$E^0_{u,v}=\{ {\rm chains~on}~\sigma [\alpha ]\in C_p(Arc):v=-{\rm rank}(H_1(F_\alpha))~{\rm
and}~u=p-v\} ,$$
where $\partial _1:E_{u,v}^0\to E_{u-1,v}^0$ and the differential of the $E^0$ term is $\partial _2:E_{u,v}^0\to
E_{u,v-1}^0$. 

\vskip .1in

\noindent
On the other hand, it follows from [20] that the $\partial _1$-homology of 
the top row itself agrees with the homology of Riemann's moduli
space of a bordered surface $F$ with one distinguished point on each boundary component (i.e., the space of hyperbolic structures on
a surface with distinguished points which is homeomorphic to $F$ modulo orientation-preserving homeomorphisms pushing foward structure and
distinguished points).  The homology of Riemann's moduli space and the arc complex are thus related by this spectral sequence, which may also
be of utility in their calculations.

\vskip .1in

\noindent One principal application of the sphericity 
results (indeed, the original impetus for
the sphericity question) is the relation with the ``combinatorial 
compactification'' of Riemann's moduli space in [19], which we next recall
(cf. also [1]). 
For any punctured surface $F_g^s$ without boundary, let us choose among the punctures
of $F_g^s$ a distinguished one.  We may then consider the arc complex comprised of all arc families
based at this distinguished puncture in the natural way.  The subspace corresponding to
quasi filling arc families was shown [17] to be real analytically identified with the usual Riemann
moduli space
$M_g^s$ of
$F_g^s$ (where the mapping class group elements are required to fix only the distinguished puncture), so the arc
complex itself forms a combinatorial compactification of $M_g^s$ 
(with no decoration but with a distinguished puncture).  As explained in [19], 
arguments analogous to the proofs of Lemma~6-7 together with sphericity
of arc complexes for polygons and multiply punctured polygons implies:

\vskip .2in

\noindent {\bf Corollary 17}~[19]~\it The combinatorial compactification of Riemann's moduli space
of a punctured surface $F_g^s$ is an orbifold in the case $g=0$ of multiply punctured spheres.\rm

\vskip .2in

\noindent 
[19] and the survey [21] contain further information about the local groups in the orbifold structure at orbifold points as well as
the geometry and combinatorics of the combinatorial compactification in general. 

\vskip .1in

\noindent Finally, it is natural to wonder what are the four special 
manifolds arising as arc complexes of the type 1 surfaces.
Athanase Papadopoulos has also asked what types of manifolds or orbifolds 
might arise for quotients of $Arc'(F)$ by
subgroups of the mapping class groups, and this has turned out to be
an excellent question.

\vskip .3in

\noindent{\bf Appendix A.  Low dimensional examples}

\vskip .2in

\noindent We explicitly describe here various low-dimensional examples of arc complexes and directly prove sphericity. 
Together with Lemma 8, these examples cover all arc complexes of dimension at most
four.

\vskip .1in

\noindent {\bf Example 1}~[Polygons]~If $g=s=0$ and $r=1$, then $F$ is a polygon with
$\delta _1$ vertices, $N=\delta _1-4$, $PMC(F)$ is trivial, and we are in the setting of the Classical Fact.  When
$\delta _1=4$ with
$F$ a quadrilateral, we have $N=0$, and there are exactly two chords, which cannot be simultaneously disjointly
embedded.  Thus,
$Arc(F)=Arc'(F)$ consists of two vertices, namely, $Arc(F)$ is the 0-sphere.  (The reader may likewise directly
investigate the case
$N=1$ where
$F$ is a pentagon, and the usual ``pentagon relation'' shows that $Arc(F)$ itself is also a pentagon.)  

\leftskip=0ex

\vskip .2in

\noindent {\bf Example 2}~~~[Once-punctured polygons]~If $g=0$, $r=s=1$, and $\delta _1<2$, then $Arc(F)=\emptyset$.  If $r=2$, then $F$ is a
once-punctured bigon.  Any arc in $F$ with distinct endpoints must be inessential, and any
arc whose endpoints coincide must be separating.  There are two
complementary components to $a$, one of which has one cusp in its
frontier, and the other of which has two cusps.  The component of the former type must
contain the puncture and of the latter type must be a triangle.  It follows that
$Arc(F)=Arc'(F)$ is in this case again a zero-dimensional sphere.

\vskip .2in\leftskip=0ex

\noindent {\bf Example 3}~[Annuli]~If $g=s=0$, $r=2$, and $\delta _1=\delta _2=1$, then $F$ is an
annulus with one distinguished point on each boundary component.  It is elementary that
any essential arc must have distinct endpoints and its isotopy class is classified by an integral ``twisting
number'', namely, the number of times it twists around the core of the annulus;
furthermore, two non-parallel essential arcs can be disjointly embedded in $F$ if and only
if their twisting numbers differ by one.  Thus, $Arc'(F)$ is in this case isomorphic to the
real line, the mapping class group $PMC(F)$ is infinite cyclic and generated by the Dehn twist along
the core of the annulus, which acts by translation by one on the real line, and the
quotient is a circle. 
\vskip .2in

\noindent {\bf Example 4}~[Twice-punctured polygons]~If $g=0$, $s=2$, and $r=1=\delta _1$, then $F$ is a
twice-punctured monogon.  If $a$ is an essential arc in $F$, then its one-cusped component
contains a single puncture, which may be either of the two punctures in $F$, and the
two-cusped component contains the other puncture of $F$.  Labeling the punctures of $F$
by 0 and 1, choose disjointly embedded arcs $a_i$ so that $a_i$ contains puncture $i$ in
its one-cusped complementary component, for $i=0,1$.  The Dehn twist $\tau$ along the
boundary generates $PMC(F)$, and the 0-skeleton of $Arc'(F)$ can
again be identified with the integers, where the vertex $\sigma (\{ \tau ^j(a_i)\} )$ is identified with
the integer $2j+i$~~for any integer $j$ and $i=0,1$.

\noindent Applying Example~2 to the
two-cusped component of $F_{\{a_i\}}$, for $i=0,1$, we find exactly two
one-simplices incident on
$\sigma (\{ a_i\} )$, namely, there are two ways to add an arc disjoint from $a_i$
corresponding to the two cusps in the two-cusped component, and these arcs are identified
with the integers $i\pm 1$.  Thus, $Arc'(F)$ is again identified with the real line,
$PMC(F)$ acts this time as translation by two, and $Arc(F)$ is again a circle.

\vskip .1in
 \noindent This complex is pictured in Figure A.1, where  the two representatives of the $PMC(F)$-orbit $[a_0]$ of $a_0$
in the illustration differ by a Dehn twist around the boundary.

\vskip .1in

\centerline{\epsffile{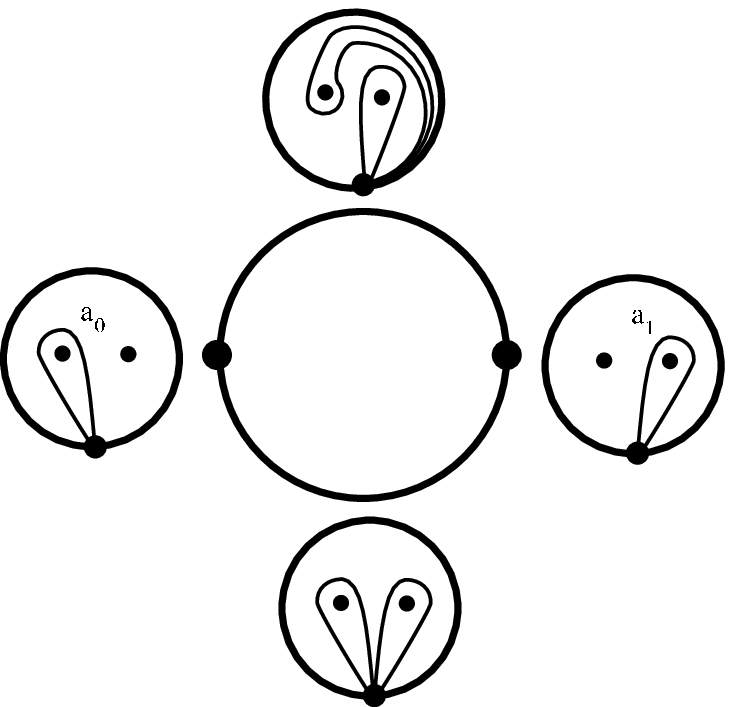}}

\centerline{{\bf Figure A.1}~The arc complex of $F_{0,1}^2$ is a circle.}
 
\vskip .1in

\vskip .2in\leftskip=0ex

\noindent {\bf Example 5} ~[Once-punctured annuli]~If $g=0$, $r=2$, and $s=1=\delta _1=\delta _2$, then $F$ is a
once-punctured annulus. To understand $Arc(F)$, first look at the subcomplex of
all arc families whose component arcs have distinct endpoints.  In a maximal such family,
the puncture lies in a bigon, and it is then easy to see that there are
exactly two 2-cells in this subcomplex of $Arc(F)$ which combine to give a torus embedded in $Arc(F)$.  The
complement of this torus retracts onto the subcomplex corresponding to arc families each of whose component arcs has
coincident endpoints.  Arguing as in Example 3, this subcomplex is isomorphic to the disjoint union of two embedded
circles (and so here is a copy of the classical Hopf fibration of $S^3$).  

\vskip .1in

\noindent More explicitly, the various cells are enumerated in Figure A.2, where subscript $\pm$ is explained as
follows.  For instance, we have depicted the arc family $\alpha_-$ in the figure, and the arc family $\alpha _+$
is obtained from $\alpha _-$ by rotating the page about the puncture in the natural way to interchange the two boundary
components.  The arcs are enumerated so that Roman letters indicate arc families where the puncture does not
lie in a monogon, and the Greek letters otherwise.

\vskip .1in

{{{\epsffile{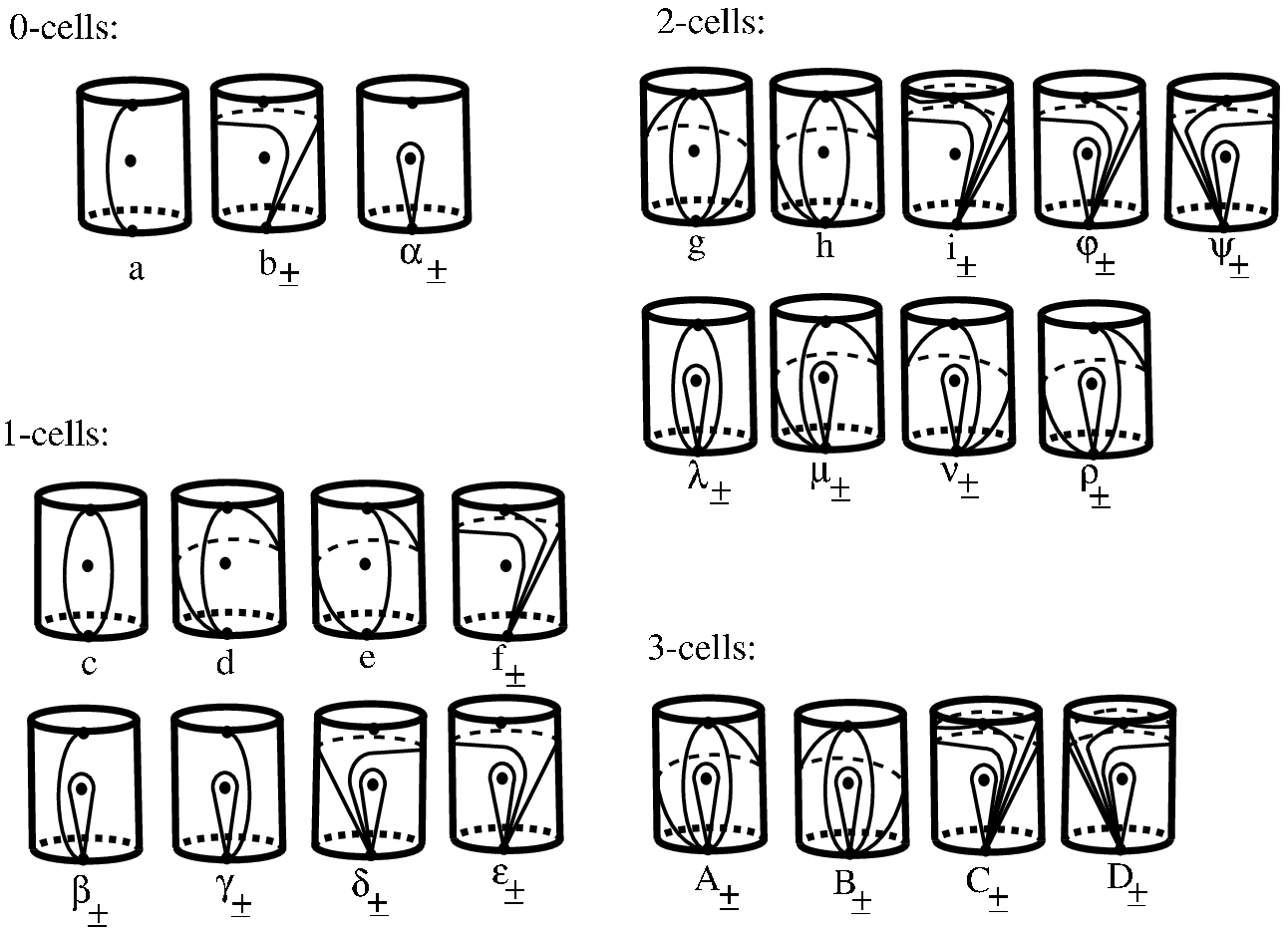}}}}

\centerline{{\bf Figure A.2}~~{The arc families in $F_{0,2}^1$.}}

\vskip .1in

~~~~~~\epsffile{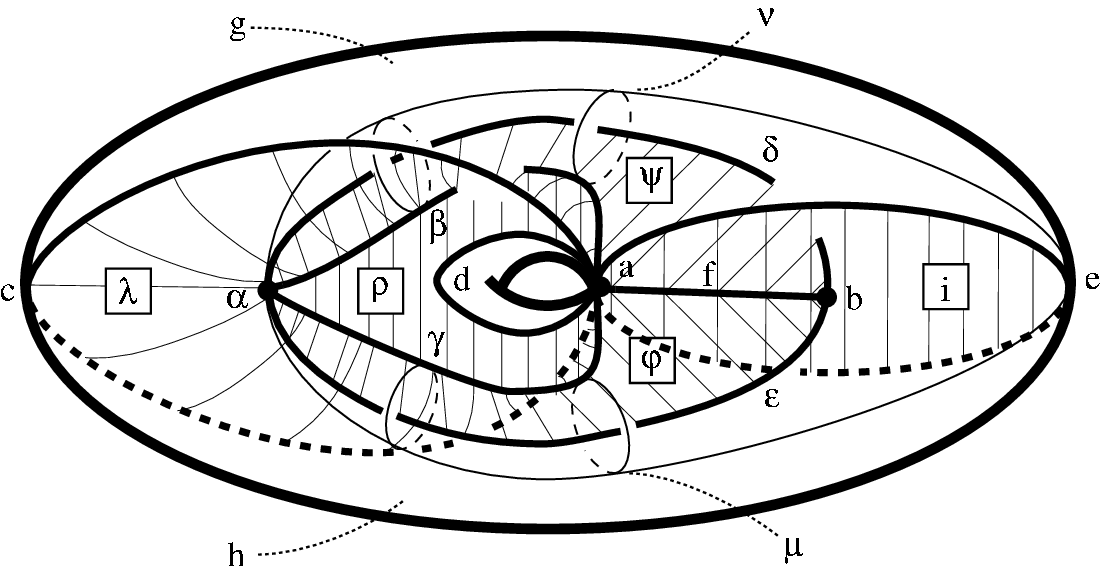}

\centerline{{\bf Figure A.3 }~~{The Hopf fibration of $Arc(F_{0,2}^1)$.}}

\vskip .2in

\noindent As illustrated in Figure A.3, the cells $a,c,d,g,h$ combine to give the torus mentioned above, the
cells with fixed subscript $\pm$ lie on one side of this torus, and we drop the subscript in the figure, where
we draw the inside of this torus.  Each 2-cell $\mu,\nu$ is a cone with vertex $\alpha$ and base $i$. 
The 3-cells $C,D$ fill these cones, and $A,B$ fill out the complement of these cones meeting along $\lambda$ and
$\rho$.  Finally, notice that the sub-complex consisting of arc families so that the puncture does not lie in a
monogon is comprised of the torus together with the two compressing disks $i_\pm\cup f_\pm\cup b_\pm$.

\vskip .2in\leftskip=0ex

\noindent {\bf Example 6}~[Torus-minus-a-disk]~ For $g=r=\delta _1=1$ and $s=0$, we shall first enumerate the cells in the
arc complex of $F=F_{1,1}^0$ and explicate the face relations.  Each arc $a$ in $F$
is non-separating,  and there is thus a unique vertex $[a]$ of $Arc(F)$.  Since
$F_{\{ a\}}$ is homeomorphic to $F_{0,(1,2)}^0$, $Arc(F_{\{
a\}})$ is completely described by Examples 1 and ~3 and consists of: four 0-cells, five
1-cells, and four 2-cells, which give rise to corresponding cells of one-higher dimension in $Arc(F)$.  It is
an exercise with $PMC(F)$ for the unflagging reader to conclude that $Arc(F)$ is comprised of:
two 1-cells corresponding to arc families
$\alpha _1, \alpha _2$, where $F_{\alpha _1}$ consists of a triangle and a surface of type
$F_{0,(1,1)}^0$, while $F_{\alpha _2}$ is a pentagon;
two 2-cells corresponding to arc families
$\beta _1, \beta _2$, where $F_{\beta _1}$ consists of a four-gon and a triangle whose frontier contains the
boundary of $F$, while $F_{\beta _2}$ consists of a triangle and a four-gon whose frontier contains the boundary
of $F$; and one 3-cell corresponding to an arc family $\gamma$, where $F_\gamma$ consists
of three triangles.  From the given descriptions, it is clear that exactly one codimension-one face of $\beta _1$
corresponds to $\alpha _1$, all three of the codimension-one faces of $\beta _2$ correspond to $\alpha _2$, and
exactly two codimension-one faces of $\gamma$ correspond to $\beta _1$.

\vskip .1in

\noindent Now consider an abstract 3-simplex with vertices $v_i$, for $0\leq i\leq 3$, and take the quotient
induced by identifying $v_0$ with $v_3$ and $v_1$ with $v_2$.  The resulting quotient space is PL-homeomorphic to
the 3-sphere, as one sees by comparing with the Hopf fibration, and has one 3-cell, two 2-cells, three 1-cells,
and two 0-cells.  Exactly one of the 1-cells has distinct endpoints, and we pass to a further quotient by
collapsing this 1-cell to a point.  The resulting quotient is again evidently PL-homeomorphic to the 3-sphere as well as
isomorphic to the complex $Arc(F)$, as desired. 

\vskip .1in

\noindent It follows from the previous paragraph that if we blow-up the
unique vertex of $Arc(F)$ to a one-simplex in the natural way, 
then we produce a space 
$\overline{Arc}(F)$ which is identified with the 3-sphere.

\vskip .1in

\noindent To close this example, we describe the cell structure of $
\overline{Arc}(F)$.
There are two types of one-simplex $\sigma [\alpha ]$ depending upon whether $\alpha$ is separating, and these two types are
represented by the weighted arc families in Figure~A.5.  In this figure, the nearby label denotes the weight depending upon the parameter
$0\leq s\leq 1$.  In the non-separating case, the two arcs $s=0$ and $s=1$ differ by a (right) Dehn twist $\tau _m$ on the meridian, and we
shall also let
$\tau _\ell$ denote the (right) Dehn twist on the longitude.  Thus, $\tau _m$ and $\tau
_\ell$ generate $PMC(F)$ and satisfy the relations $\iota=\tau _m\circ \tau _\ell \circ \tau _m=\tau _\ell\circ \tau
_m\circ \tau _\ell$, where $\iota ^4$ is represented by the (right) Dehn twist on the boundary.  In the separating case, the two arcs $s=0$ and
$s=1$ differ by the ``half Dehn twist'' $\iota ^2$.

\hskip .6in{{{\epsffile{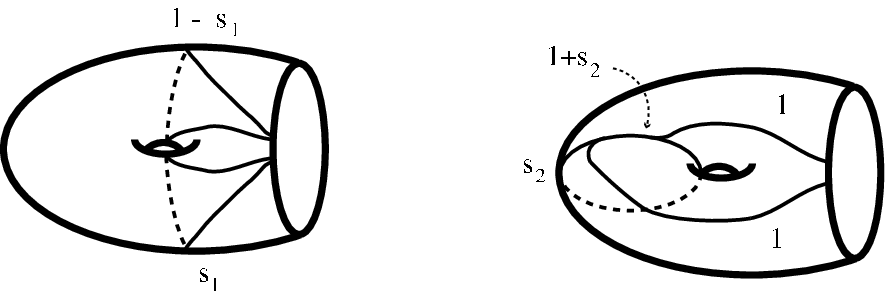}}}}

\hskip .35in{\bf A.5a}~Separating case, circle $C_1$.~~~~~~~~~~~{\bf A.5b}~Non-separating case, circle
$C_2$.

\vskip .1in

\centerline{{\bf Figure A.5}~~Two circles in $Arc(F_{1,1}^0)$.}

\vskip .1in

\noindent Thus, there are two circles $C_1,C_2\subseteq {Arc}(F)$ as illustrated in Figures A.5a and A.5b respectively.  It is
not difficult to see that $C_1,C_2$ are embedded in $Arc(F)$ and $C_1\cap C_2=\{  [a]\}$.  We may regard each $C_i$ as a CW complex with one
vertex and one edge, and there is then an induced CW complex structure on $C_1*C_2$.

\vskip .1in

\noindent We claim that $\overline{Arc}(F)$ is isomorphic to $C_1*C_2$ as CW complexes.  To see this, first observe that a point of
$C_i$ is given by a coordinate
$s_i\in S^1=[0,1]/(0\sim 1)$, for $i=1,2$, so a point of $C_1*C_2$ may be identified with the projectivization of four real
numbers $(L_1,L_2,S_1,S_2)$, where $L_1,L_2\geq 0$ and at least one of $L_1,L_2$ is non zero, and $S_i=s_iL_i\in [0,L_i]/(0\sim
L_i)$, for $i=1,2$; a parameter along the join line in the join structure of $C_1*C_2$ is given by $L_1/(L_1+L_2)$.

\vskip .1in

\noindent We may construct a
weighted arc family from this data provided $s_1\neq 0\neq s_2$ as illustrated in Figure~A.6a.  This weighted arc family projectivizes
to the unique 3-cell in $Arc(F_{1,1}^0)$.  In Figures A.6b-c are likewise illustrated the 2-cells $(s_1=0)*C_2$ and
$C_1*(s_2=0)$ which projectivize to the 2-cells
$\beta _1,\beta _2$, and the 1-cells $C_1,C_2$ in Figure~A.5 furthermore projectivize to the 1-cells $\alpha _1,\alpha _2$.  The
join parameter
${{L_1}\over{L_1+L_2}}$ provides the blow-up parameter for $\overline{Arc}(F)$, i.e., $[a]$ is replaced by the interval
$(s_1=0)*(s_2=0)$ of join parameters.

\vskip .2in

\hskip -.2in\epsffile{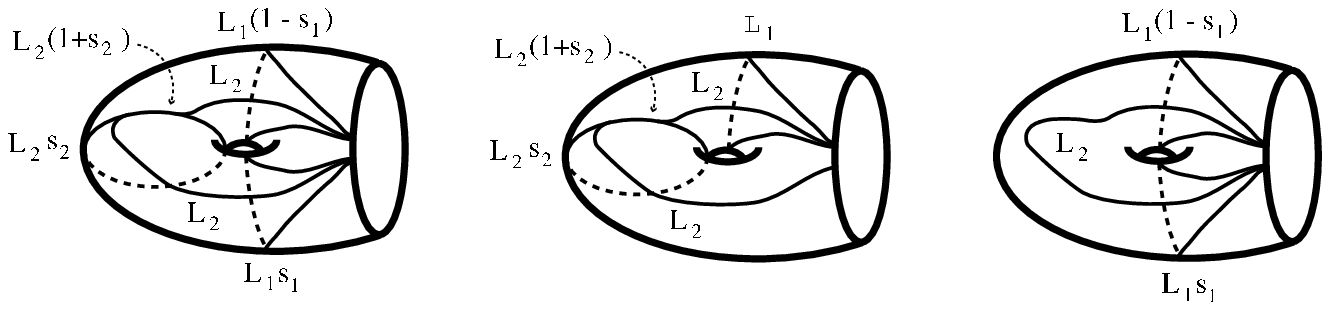}

\hskip .3in {\bf A.6a}~The 3-cell.\hskip .9in{\bf A.6b}~The 2-cell $\beta _1$.\hskip .6in{\bf A.6c}~The 2-cell $\beta _2$.

\vskip .1in

\centerline{{\bf Figure A.6}~Cell structure of the join $C_1*C_2$.}

\vskip .2in

\vskip .2in\leftskip=0ex

\noindent {\bf Example 7}~[Thrice-punctured polygons]~Label the various punctures of $F=F_{0,1}^3$ with distinct
members of the set
$S=\{ 1,2,3\}$.  If $\alpha $ is an arc family in $F$ and $a$ is a component arc of $\alpha $,
then $a$ has a corresponding one-cusped component containing some collection of punctures
labeled by a proper non-empty subset
$S(\alpha )\subseteq S$,
which we regard as the ``label'' of $a$ itself. 

\vskip .1in

\noindent More generally, define a ``tableaux'' $\tau$ labeled by $S$ to be a rooted tree embedded in
the plane where: the (not necessarily univalent) root of $\tau$ is an unlabeled vertex, and the other
vertices of $\tau$ are labeled by proper non-empty subsets of
$S$; for any $n\geq 1$, the vertices of $\tau$ at distance $n$ from the root are pairwise disjoint  subsets
of $S$; and, if a simple path in $\tau$ from the root passes consecutively through
the vertices labeled
$S_1$ and $S_2$, then $S_2$ is a proper subset of $S_1$.

\vskip .1in

\noindent Given an arc family $\alpha $ in $F$, we inductively define the corresponding
tableaux
$\tau =\tau (\alpha )$ as follows.  For the basis step, choose as root some point in the component of
$F_{\alpha }$ which contains the boundary of $F$.  For the inductive step, given a vertex of $\tau$
lying in a complementary region $R$ of $F_{\alpha }$, enumerate the component arcs
$a_0,a_1,\dots ,a_m$ of $\alpha $ in the frontier of $R$, where we assume that these
arcs occur in this order in the canonical linear ordering 
described in $\S$2 and $a_0$ separates $R$ from the root.  Each arc $a_i$ separates $R$ from another component
$R_i$ of
$F_{\alpha }$, and we adjoin to $\tau$ one vertex in each such component $R_i$ with the label $S(a_i)$ together
with a one-simplex connecting $R$ to $R_i$, for each
$i=1,2,\dots ,m$, where the one-simplices are disjointly embedded in $F$.  This completes the
inductive definition of
$\tau =\tau (\alpha )$, which evidently satisfies the conditions required of a tableaux. 
Furthermore, the isomorphism class of
$\tau$ as a labeled rooted tree in the plane is obviously well-defined independent of the
representative
$\alpha
$ of the $PMC(F)$-orbit $[\alpha ]$.

\vskip .1in

\noindent It follows immediately from the topological classification of surfaces that
$PMC(F)$-orbits of cells in $Arc(F)$ are in one-to-one
correspondence with isomorphism class of tableaux labeled by $S$.
Furthermore, since the edges of $\tau (\alpha )$ are in one-to-one correspondence with the component
arcs of
$\alpha $, a point in $Int(\sigma (\alpha ))$ is uniquely
determined by a projective positive weight on the edges of $\tau(\alpha )$.  It follows that $Arc(F)$
itself is identified with the collection of all such projective weightings on all isomorphism classes of
tableaux labeled by $S$.

\vskip .1in

\noindent Let us adopt the convention that given an ordered pair $ij$, where $i,j\in\{
1,2,3\}$, we shall let $k=k(i,j)=\{ 1,2,3\} -\{ i,j\}$, so $k$ actually depends only upon the unordered pair
$i,j$.  The various tableaux for $F_{0,1}^3$ are enumerated, labeled, and indexed in Figure~A.7, where in each case,
$ij$ varies over all ordered pairs of distinct members of $\{ 1,2,3\}$, $k=1,2,3$, and the bullet represents the
root.   In this notation and
letting
$\partial$ denote the boundary mapping in $Arc=Arc(F_{0,1}^3)$, one may directly compute incidences of cells summarized as follows.

\settabs 5\columns

\vskip .2in

\+&$\partial C_{ij}=A_i-A_j$,&&$\partial G_{ij}=C_{ij}-D_{ji}+D_{ij}$,\cr
\+&$\partial D_{ij}=A_j-B_k$,&&$\partial H_{ij}=D_{ij}-C_{jk}+F_k$,\cr
\+&$\partial E_{k}=A_k-B_k$,&&$\partial I_{ij}=D_{ij}-E_k+C_{kj}$,\cr
\+&$\partial F_{k}=B_k-A_k$,&&$\partial J_{ij}=C_{ij}-C_{ik}+C_{jk}$,\cr
\+&$\partial K_{ij}=G_{ij}-H_{ij}+H_{ji}-J_{ij}$,&&$\partial L_{ij}=I_{ij}-G_{ij}+J_{ki}-I_{ji}$.\cr

\vskip .2in

\hskip 1.in {{{\epsffile{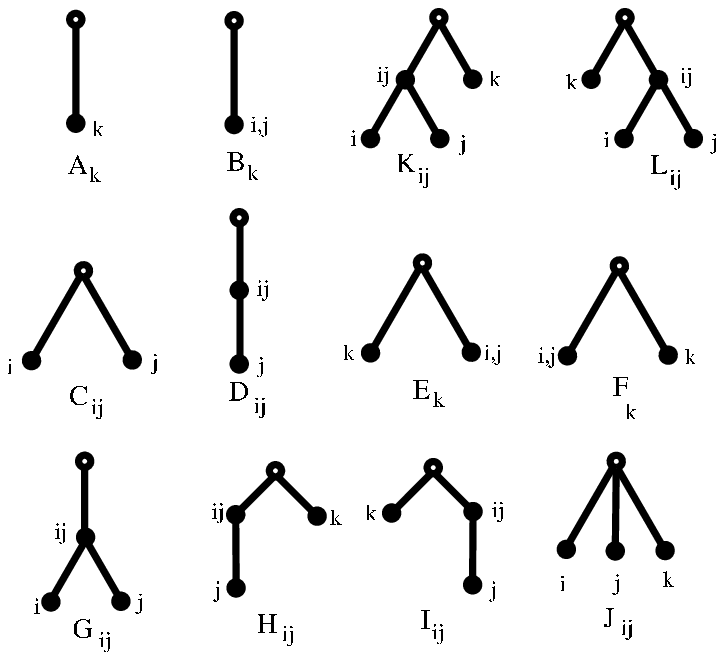}}}}

\centerline{{\bf Figure A.7}~~{The tableaux for $F_{0,1}^3$.}}

\vskip .2in

\noindent We may symmetrize and define sub complexes ${X}_{ k}=X_{ij}\cup X_{ji}$, for $X=K,L$, and furthermore set
${M}_{ k}={ K}_{k }\cup {L}_{k}$, for $k=1,2,3$.  Inspection of the incidences of cells shows that each of ${K}_k$ and
${L}_k$ is a 3-dimensional ball embedded in $Arc$, as illustrated in Figures~A.8a and 8b, respectively,
with $K_{ij}$ on the top in part a) and $L_{ji}$ on the top in part b).  Gluing together $L_k$ and $K_k$ along their common
faces $G_{ij}, G_{ji}$, we discover that
${M}_{k}$ is almost a 3-dimensional ball embedded in
$Arc$ except that two points in its boundary are identified to the single point $A_k$ in $Arc$, as illustrated in Figure~A.8c.

\vskip .2in

~{{{\epsffile{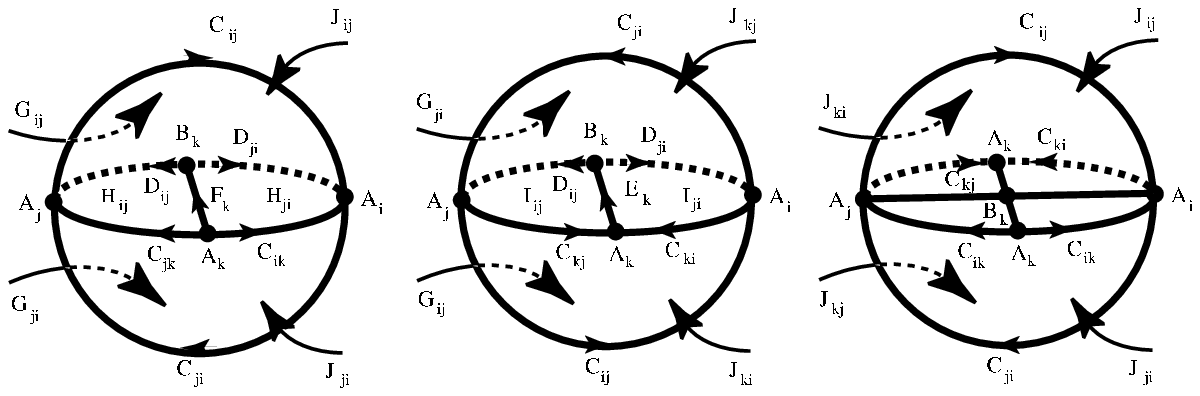}}}}

\vskip .1in

\noindent \hskip .65in Figure~A.8a ${K}_{k}$.\hskip .6in Figure A.8b ${ L}_{k}$.\hskip .6in Figure A.8c ${
M}_{k}$.

\vskip .1in

\centerline{{\bf Figure A.8}~~{${ M}_{k}$ and the balls ${ K}_{k},{ L}_{k}$.}}

\vskip .2in

\noindent Each ${ M}_k$, for $k=1,2,3$, has its boundary entirely contained in the sub complex
${J}$ of $Arc$ spanned by $$\{ A_k: k=1,2,3\}\cup\{ C_{ij},J_{ij}: i,j\in\{ 1,2,3\}~{\rm are~distinct}\}.$$  In order to understand
${J}$, we again symmetrize and define $J_{k}=J_{ij}\cup J_{ji}$, so each $J_{k}$ is isomorphic to a cone, as
illustrated in Figure A.9; 
we shall refer to the 1-dimensional simplices $C_{ik},C_{jk}$ as the ``generators'' and to
$C_{ij},C_{ji}$ as the ``lips'' of
$J_{k}$.  The one-skeleton of ${ J}$ plus the cone $J_{k}$ is illustrated in Figure A.9.
Imagine taking $k=3$ in Figure A.9 and adjoining the cone $J_2$ so that  
the generator $C_{12}$ of $J_{2}$ is attached to the lip $C_{12}$ of $J_{3}$, and the generator $C_{13}$ of
$J_{3}$ is attached to the lip $C_{13}$ of $J_{2}$.  Finally, ${ J}$ itself is produced by symmetrically attaching
$J_{1}$ to $J_{2}\cup J_{3}$ in this lip-to-generator fashion.

\hskip 1.3in{{{\epsffile{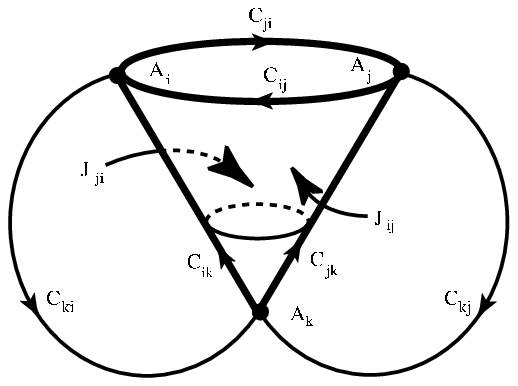}}}}

\vskip .2in

\centerline{{\bf Figure A.9}~~{The cone ${ J}_{k}$ in the one-skeleton of ${ J}$.}}

\vskip .2in

\noindent In order to finally recognize the 3-dimensional sphere, it is best to take a regular neighborhood of ${ J}$ in
$Arc$, whose complement is a disjoint union of three 3-dimensional balls.  Each 3-dimensional ball is naturally identified with
the standard ``truncated '' 3-simplex, where a polyhedral neighborhood of the 1-skeleton of the standard 3-simplex has been excised. 
These truncated simplices are identified pairwise along pairs of faces to produce the 3-dimensional sphere in the natural way.
This completes Example~7 and our serial explicit treatment of low-dimensional examples.

\vskip .2in

\noindent {\bf Proposition A}~\it The examples of this appendix establish
sphericity of all arc complexes of dimension at most four.\rm

\vskip .2in

\noindent{\bf Proof}~~By Lemma~8, we may assume that $r=\Delta$.  Enumerating the solutions to
$6g-7+4r+2s\leq 4$, one finds precisely the examples of this appendix.~~~\hfill{\it q.e.d.}

\vskip .3in

\noindent{\bf Appendix B.  Pairs of pants}

\vskip .2in

\noindent It is the purpose of this section to prove sphericity of the arc complex for the ``generalized pairs of pants'' $F_{0,r}^s$, where
$r+s=3$.  

\vskip .1in

\noindent To this end, fix a pair of pants $P=F_{0,3}^0=F_{0,(1,1,1)}^0$ and consider a non-empty family $\alpha$ of disjointly embedded
essential arcs with distinct endpoints in the boundary of $P$.  (Components of $\alpha$ may be parallel to one another, so $\alpha$ is not an
arc family in the sense of $\S$1, and furthermore, there are twice the cardinality of $\alpha$ many points in $(\cup\alpha )\cap \partial
P$.)  Label the components of the boundary $\partial P$ of $P$ as
$\partial _1,\partial _2,\partial _3$, and define the {\it intersection number} $m_i\geq 0$ to be the cardinality of $\partial _i\cap
(\cup\alpha )$, i.e., the number of endpoints of components of $\alpha$ lying in $\partial _i$, for $i=1,2,3$.  Provided $\alpha
\neq\emptyset$, at least one of $m_1, m_2, m_3$ is non-zero,
and furthermore, 
$m_1+m_2+m_3$ must be even since each component arc in $\alpha$ has two endpoints.
More generally, any weighted arc family in $P$ has corresponding total weights on $\partial _i$ also denoted $m_i\in {\bf R}$,
for $i=1,2,3$.

\vskip .2in

~~~~\epsffile{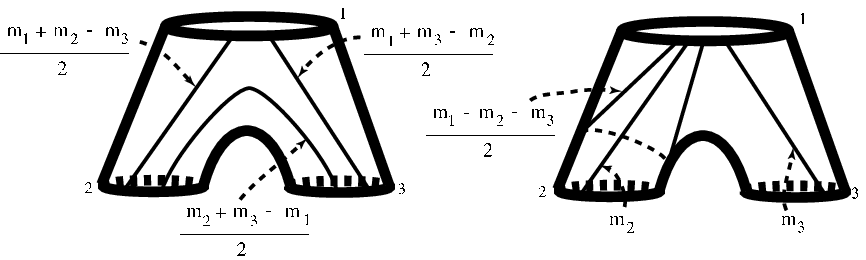}

\hskip .4in{\bf B.1a}~Triangle inequality\hskip .5in{\bf B.1b}~The case $m_1>m_2+m_3$

\vskip .1in

\centerline{{\bf Figure B.1}~Constructing arcs in pants}

\vskip .1in

\noindent {\bf Dehn-Thurston Lemma}~\it Isotopy classes of (non-empty) families of disjointly embedded essential arcs
in $P$ are uniquely determined by the triple $(m_1,m_2,m_3)$ of non-negative integral intersection numbers, which are subject only to the
constraints that
$m_1+m_2+m_3$ is positive and even.  More generally, isotopy classes of weighted arc families in $P$ are uniquely determined by the
analogous real triple $(m_1,m_2,m_3)$ of nonnegative total weights, which are subject to the unique constraint that $m_1+m_2+m_3$ is
positive.\rm

\vskip.2in

\centerline{\epsffile{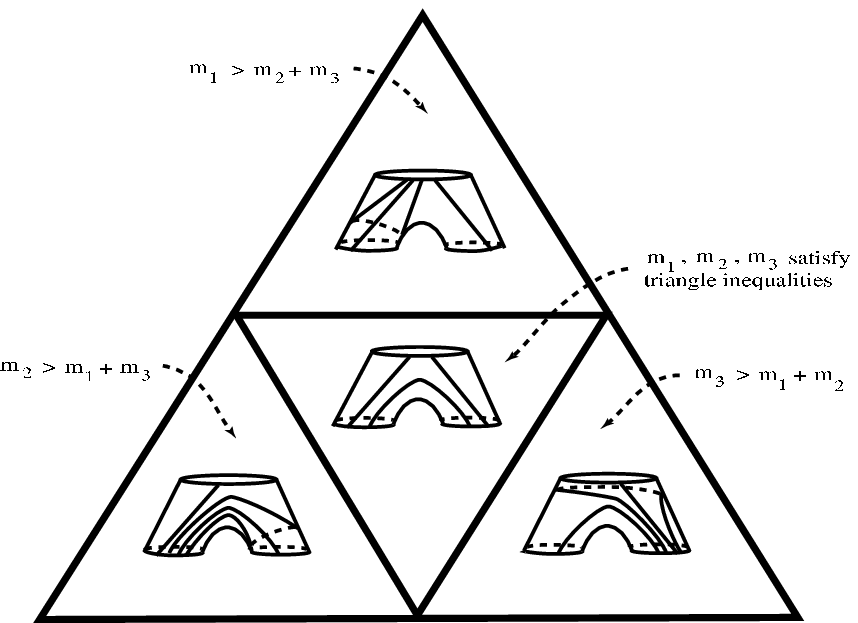}}

\vskip .1in

\centerline{{\bf Figure B.2}~Families of arcs in pants}

\vskip .1in

\noindent {\bf Proof}~~The explicit construction of a family of arcs realizing a putative triple of intersection numbers
is illustrated in Figure~B.1 in the two representative cases that $m_1,m_2,m_3$ satisfy all three possible (weak) triangle inequalities
in Figure~B.1a
or perhaps $m_1\geq m_2+m_3$ in Figure~B.1b; notice that the weights in Figure~B.1 are integral by the parity condition on $m_1+m_2+m_3$.  The
other cases are similar, and the projectivization of the positive orthant in
$(m_1,m_2,m_3)$-space is illustrated in Figure~B.2.  Elementary topological considerations show that any family of arcs in $P$ is isotopic
to a unique such family keeping endpoints of arcs in the boundary of $P$, completing the proof.~~~~\hfill{\it q.e.d.}

\vskip .2in

\noindent In order to refine The Dehn-Thurston Lemma and keep track of twisting around the boundary, we shall introduce in each
component
$\partial _i$ of the boundary an arc $w_i\subseteq \partial _i$ called a ``window'', for $i=1,2,3$.  We require arcs in a family to
have their endpoints in the windows, and we consider isotopy classes of such families of arcs rel windows and shall call these
``windowed'' isotopy classes of ``windowed'' families of arcs.  
(Collapsing each window to a point gives a surface with one cusp on each boundary component, so a windowed arc family gives
rise  to an arc family in the usual sense of $\S$1.)  We seek the analogue of the Dehn-Thurston Lemma for windowed isotopy classes.

\vskip .1in

\noindent To this end, there are two conventions to be made:

\vskip .1in

\leftskip .3in

\noindent 1) when an essential arc in 
$P$ connects a boundary component to itself, i.e., when it is a loop, then it passes around a specified leg of
$P$ as illustrated in Figure B.3a-b, i.e., it contains a particular boundary component in its one-cusped component;

\vskip .1in 

\noindent 2) when an arc is added to the two-cusped complementary region of a loop in $P$, then it either
follows or precedes the loop as illustrated in Figure B.3c-d.

\leftskip=0ex

\vskip .2in 

\noindent For instance in Figure~B.2, the conventions are: 1) around the right leg
for loops; and 2) on the front of the surface.

\vskip .05in

\hskip .1in{{{\epsffile{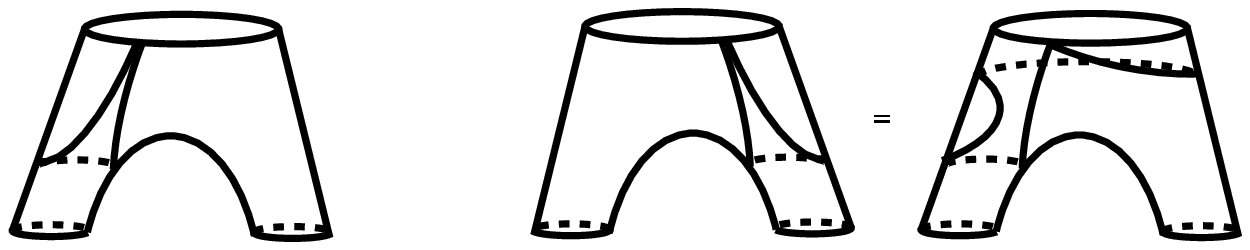}}}}

\hskip .1in{\bf Figure B.3a} Around the right leg.\hskip .5in{\bf Figure
B.3b} Around the left leg.

\vskip .05in

\hskip .1in{{{\epsffile{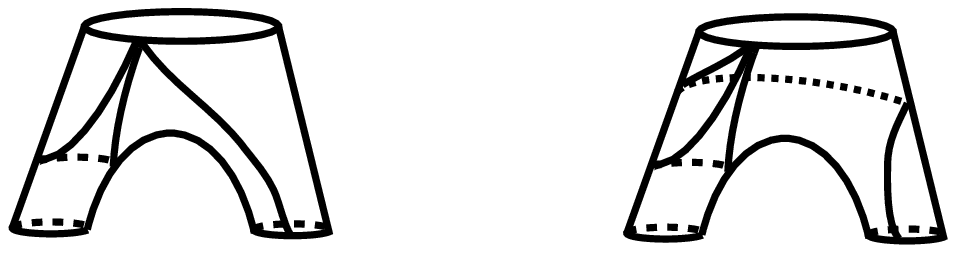}}}}

\hskip .1in{\bf Figure B.3c} Arc follows loop.\hskip .5in~~~~~~{\bf Figure
B.3d}~Arc precedes loop.

\vskip .05in

\centerline{{\bf Figure B.3}~~{Twisting conventions.}}

\vskip .1in

\noindent Upon making such choices of convention, we may associate a twisting number $t_i\in {\bf Z}$ to each family of arcs and
$t_i\in{\bf R}$ to each weighted arc family as follows.  Choose a regular neighborhood of $\partial P$ and consider the sub-pair of
pants
$P_1\subseteq P$ complementary to this regular neighborhood.  Given a weighted arc family $\alpha$ in $P$, by the Dehn-Thurston Lemma,
we may perform a windowed isotopy in $P$ supported on a neighborhood of $P_1$ to arrange that $\alpha\cap P_1$ agrees with a conventional
windowed arc family in $P_1$ (where the window in $\partial P_1$ arises from that in $\partial P$ in the natural way via a framing of the
normal bundle to $\partial P$ in $P$).  

\vskip .1in
 \noindent For such a
representative
$\alpha$, we finally consider its restriction to each annular neighborhood $A_i$ of $\partial _i$.  As in Example~3 in Appendix~A,
there is an associated real quantity $t_i$, called a {\it twisting number} defined explicitly as follows: choose another arc $a$ 
whose endpoints are {\sl not} in the windows (and again such an arc is essentially uniquely determined up to
windowed isotopy from a framing of the normal bundle to
$\partial P$ in $P$ in the natural way); orient $a$ and each component arc of $(\cup \alpha )
\cap A_i$ from
$\partial P_1$ to $\partial P$, and let $t_i$ be the signed (weighted) intersection number of $a$ with the (weighted)
arc family $(\cup \alpha)\cap A_i$, for $i=1,2,3$.

\vskip .1in

\centerline{\epsffile{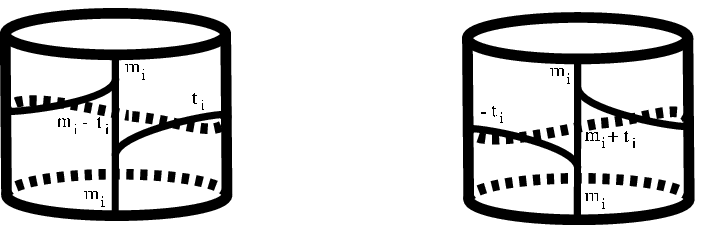}}

\vskip .1in

\hskip .6in{\bf B.4a}~Right twisting for $t_i\geq 0$\hskip .3in{\bf B.4b}~Left twisting for $t_i\leq 0$

\vskip .1in

\centerline{{\bf Figure B.4}~Windowed arc families in the annulus}

\vskip .1in

\noindent As illustrated in Figure~B.4, where the weight next to an edge
indicates the width of a corresponding band in the weighted arc family, all possible real twisting numbers $-m_i\leq t_i\leq
m_i$ arise provided
$m_i\neq 0$.  By performing Dehn twists along the core of the annulus, it likewise follows that every real twisting number $t_i$
occurs provided $m_i\neq 0$.  Again, elementary topological considerations show that each windowed isotopy class of a windowed family of
arcs is uniquely determined by its invariants:

\vskip .2in

\noindent {\bf Lemma B}~\it Windowed isotopy classes of (non-empty) families of disjointly embedded essential windowed arcs
in $P$ are uniquely determined by the triple $(m_i,t_i)$ of pairs of nonnegative integral intersection number
and integral twisting number, which are subject
only to the constraints that
$m_1+m_2+m_3$ is positive and even, and $\forall i=1,2,3(m_i=0\Rightarrow t_i=0)$.  More generally, classes of windowed weighted
arc families in
$P$ are uniquely determined by the analogous triples $(m_i,t_i)$ of nonnegative real intersection number and real twisting number, which are
subject only to the constraints that
$m_1+m_2+m_3$ is positive, and $\forall i=1,2,3(m_i=0\Rightarrow t_i=0)$.\rm

\vskip .1in

\noindent {\bf Remark}~In the related setting [5,18,24,27] of measured 
foliations in $P$, the constraint $m_i=0\Rightarrow t_i=0$ is
relaxed to the constraint $m_i=0\Rightarrow t_i\geq 0$, where the coordinate $m_i=0$ and $t_i=|t_i|> 0$ corresponds to the
class of a foliated annulus of width $t_i$ whose leaves are parallel to $\partial _i$.  Thus, in the topology of projective measured
foliations, extensive twisting to the right or left about $\partial _i$ approaches the curve parallel to $\partial _i$.  More generally,
choosing a collection of curves in an arbitrary surface
$F$ which decompose
$F$ into ``generalized pairs of pants'', we may likewise assign an intersection and twisting number to each curve in the collection to give a
parameterization of the space of all classes of measured foliations in $F$; see [18] or [24] for more information, which also give explicit
formulas for the action of Dehn twist generators for the mapping class group on these global coordinates for the space of measured foliations
in any surface.  

\vskip .1in

\noindent The pure mapping class group $PMC(P)$ of $P$
is generated by the three commuting Dehn twists along the boundary components.  The
action of the Dehn twist along the $i^{\rm th}$ boundary component is
given in coordinates by
$t_i\mapsto t_i\pm m_i$ with the other coordinates unaffected.
The resulting ${\bf Z}^3$-action has a fundamental domain
described by $|t_i|\leq m_i/2$ as one can see directly by referring again to Figure~B.4.

\vskip .1in

\noindent More explicitly, there is a canonical projection from the pre-arc complex $Arc'(P)$ of $P$ to the
two-simplex $\Sigma ^2$ of projective classes of intersection numbers which is illustrated in Figure~4; the fiber over a point in
the interior $\Sigma ^2$ is identified with ${\bf R}^3$, over a one-simplex in the frontier of $\Sigma ^2$ with ${\bf R}^2$, and over a vertex
in the frontier with ${\bf R}$.  
This gives a ``toric
bundle'' structure to $Arc(P)=Arc'(P)/PMC(P)$ in the sense that the fiber over an interior point of $\Sigma ^2$ is identified with a 3-torus
$S^1\times S^1\times S^1$, over a one-simplex in the frontier with a 2-torus $S^1\times S^1$, and over a vertex in the frontier with a circle
$S^1$.  

\vskip .1in

\noindent The toric bundle $Arc(P)\to \Sigma ^2$ thus has a well-defined section $t_1,t_2,t_3$ determined by the twisting conventions.
Changing the twisting conventions affects the section by an overall rotation of the corresponding circles, where the rotation depends only
upon $m_1,m_2,m_3$.  In fact, the toric bundle structure is induced from a join structure on $Arc(P)$, as follows.

\vskip .1in

\noindent In
the $m_i,t_i$ plane, consider the region where $m_i> 0$ and $|t_i|\leq m_i/2$; identify the two rays in the
frontier of this region to get a cone minus its apex; the projectivization of this cone is an abstract circle
$S^1_i$ containing a coordinate, whose deprojectivization gives a pair $(m_i,t_i)$ of (intersection, twisting)
numbers, for $i=1,2,3$.

\vskip .1in

\noindent {\bf Proposition B}\it ~$Arc(P)=Arc(F_{0,3}^0)$ is homeomorphic to $S_1^1 * S^1_2 *
S_3^1$ and hence has a natural toric bundle structure.  Likewise, $Arc(F_{0,2}^1)\approx S^1*S^1$ and $Arc(F_{0,1}^2)\approx S^1$ have natural
induced structures.\rm

\vskip .1in

\noindent {\bf Proof}~A point of the deprojectivized arc complex
of $P$ is uniquely determined from Lemma~B by a tuple of deprojectivized coordinates $(m_i,t_i)$, where $m_1+m_2+m_3\neq 0$, as
required for $Arc(F_{0,3}^0)$.  The other surfaces
$F_{0,r}^s$, for $r+s=3$ then follow from Lemma~2.~~~~~\hfill{\it q.e.d.}

\hskip 1.2in\epsffile{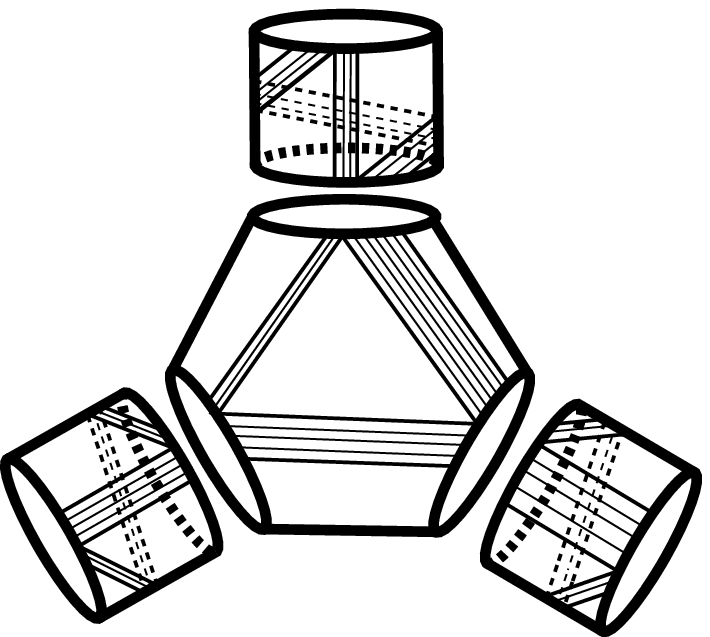}

\vskip .1in

\centerline{{\bf Figure B.5}~~Explicit construction of arcs in pants.}

\vskip .1in

\noindent  Since it is so simple, it is worth explicitly describing the direct construction of weighted arc families in $P$ 
from the $(m_i,t_i)$ coordinates, for $i=1,2,3$.  For each $m_i\neq 0$, construct a weighted arc family in the standard annulus
realizing $m_i,t_i$ as intersection and twisting number, for $i=1,2,3$, as in Figure~B.4; also, construct the conventional weighted arc
family in $P$ as in Figures~B.2-3 realizing the intersection numbers $m_1,m_2,m_3$.  Attach these annuli to $P$ in the natural way
gluing band representatives of weighted arcs to band representatives of weighted arcs without twisting to produce the required
weighted arc family in $P$ as illustrated in Figure~B.5.

\vskip .3in

\vfill\eject

\noindent {\bf Bibliography}

\vskip .2in

\noindent [1]~Brian H.
Bowditch and D. B. A. Epstein,
``Natural triangulations associated to a surface'',
{\it Top.} {\bf 27} (1988), 91-117.

\vskip .1in

\noindent [2]~Moira Chas and Dennis P. Sullivan, ``String topology'', Preprint math.GT/9911159,
to appear {\it Ann. Math.}.

\vskip .1in

\noindent [3]~---, 
``Closed string operators in topology leading to Lie bialgebras and 
higher string algebra'',  The legacy of Niels Henrik Abel,  
Springer, Berlin (2004), 771-784, Preprint math.GT/0212358.

\vskip .1in

\noindent [4]~Ralph L.\ Cohen and John D.S.\ Jones
{\it  A homotopy theoretic realization of string topology}
Preprint math.GT/0107187.

\vskip .1in

\noindent [5] Albert Fathi, Francois Laudenbach, Valentin Poenaru, 
``Travaux de Thurston sur les surfaces'', {\it Asterisque} {\bf 66-67},
Soc. Math. de France (1979).

\vskip .1in

\noindent [6] John L. Harer, 
``The virtual cohomological dimension of the mapping class group of an 
orientable surface'',  {\it Invent. Math.}  {\bf 84}  (1986),  157-176. 

\vskip .1in

\noindent [7] ---,
``Stability of the homology of the mapping class groups of orientable 
surfaces'',  {\it Ann. Math.}  {\bf  121}  (1985), 215-249.

\vskip .1in
                                                                                
\noindent [8] John L. Harer and Don Zagier,
``The Euler characteristic of the moduli space of curves'',
{\it Invent. Math.} {\bf 85} (1986), 457-485.

\vskip .1in
                                                                                
\noindent [9] William J. Harvey,
``Boundary structure of the modular group'',
Riemann surfaces and related topics: Proceedings of the 1978 Stony Brook
Conference,  {\it Ann. of Math. Stud.} {\bf 97},
Princeton Univ. Press, Princeton, N.J., (1981), 245-251.

\vskip .1in

\noindent [10] John H. Hubbard and Howard Masur,
``Quadratic differentials and foliations'',  {\it Acta Math.}  {\bf 142}
(1979), 221-274.

\vskip .1in
                                                                                
\noindent [11] Kyoshi Igusa,
``Combinatorial Miller-Morita-Mumford classes and Witten cycles'',
{\it  Algebr. Geom. Topol.}  {\bf 4}  (2004), 473-520.
                                                                                
\vskip .1in
                                                                                
\noindent [12] Nikolai Ivanov,
``Stabilization of the homology of Teichm\"uller modular groups''
(Russian)  {\it Algebra i Analiz} {\bf  1}  (1989), 110-126;
translation in  {\it Leningrad Math. J.} {\bf 1}  (1990), 675-691.
                                                                                
\vskip .1in
                                                                                
\noindent [13] ---,
``On the homology stability for Teichm\"uller modular groups: closed
surfaces and twisted coefficients'', Mapping class groups and moduli
spaces of Riemann surfaces (G\"ottingen, 1991/Seattle, WA, 1991),
{\it Contemp. Math.} {\bf 150}, Amer. Math. Soc., Providence, RI, 1993, 149-194.
                                                                                
\vskip .1in

\noindent [14]~Ralph L. Kaufmann, Muriel Livernet, R. C. Penner, 
``Arc operads and arc algebras'', 
{\it Geom. Topol.} {\bf 7} (2003), 511-568,
Preprint math.GT/0209132.

\vskip .1in
                                                                                
\noindent [15] Maxim Kontsevich,
``Intersection theory on the moduli space of curves and the
matrix Airy function'',  {\it Comm. Math. Phys.}  {\bf  147}  (1992), 1-23.

\vskip .1in
                                                                                
\noindent [16] Gabriel Mondello,
``Combinatorial classes
on $\overline{M}\sb {g,n}$ are tautological'',  {\it Int. Math. Res. Not.}
(2004), 2329-2390.

\vskip .1in

\noindent [17]~R. C. Penner, ``The decorated Teichm\"uller space of  punctured surfaces", 
{\it Comm. Math.  Phys.}  {\bf 113}   (1987),  299-339.
\vskip .1in

\noindent [18]~---, ``The action of the mapping class group on curves and arcs in surfaces'', thesis,
Mass. Inst. of Tech. (1982).

\vskip .1in

\noindent  [19]~---,``The simplicial compactification of Riemann's moduli space'', 
Proceedings of the 37th Taniguchi Symposium, World Scientific (1996), 237-252.

\vskip .1in

\noindent [20]~---, ``Decorated Teichm\"uller theory of bordered surfaces'', {\it Comm. Anal.  Geom.} {\bf 12} (2004), 793-820,
math.GT/0210326.

\vskip .1in

\noindent [21]~---, ``Cell decomposition and compactification of Riemann's moduli space in decorated Teichm\"uller theory, 
Woods Hole volume, World Scientific Press (2003).

\vskip .1in
                                                                                
\noindent [22] ---,
``Perturbative series and the moduli space of Riemann surfaces'',
{\it  J. Diff. Geom.}  {\bf 27}  (1988), 35-53.
                                                                                
\vskip .1in
                                                                                
\noindent [23] ---,
``Weil-Petersson volumes'',
{\it  J. Diff. Geom.} {\bf  35}  (1992), 559-608.
                                                                                
\vskip .1in

\noindent [24]~---~with John L. Harer, ``Combinatorics of Train Tracks'', Annals of Math Studies {\bf 125}, {Princeton
Univ. Press} (1992).

\vskip .1in

\noindent [25]~---~and Michael Waterman, ``Spaces of RNA secondary structures",
{\it Adv. Math.}  {\bf 101}  (1993), 31-49.

\vskip .1in
                                                                                
\noindent [26] Kurt Strebel, Quadratic Differentials, {\it Ergebnisse
der Math. und ihrer Grenzgebiete}, Springer-Verlag, Berlin (1984).
                                                                                
\vskip .1in

\noindent [27]~Wlliam~P.~Thurston, ``On the geometry and dynamics of diffeomorphisms of surfaces'',
{\sl Bull. Amer. Math. Soc.}, {\bf 19} (1988) 417--431.

\vskip .1in

\noindent [28]~Alexander 
A. Voronov {``Notes on universal algebra''}, Preprint math.QA/00111009.

\vskip .1in

\noindent [29]~John Milnor, ``Lectures on the h-Cobordism Theorem'' (with L. Siebenmann and J. Sondow), Princeton 
U. Press, 1965.

\vskip .1in

\noindent [30]~Armand Borel, ``Seminar on Transformation Groups'', Annals of Math Studies
{\bf 46}, Princeton University Press, 1960.

\vskip .1in

\noindent [31]~Marhall M. Cohen, ``Simplicial structures and transverse cellularity'', {\it Ann. Math.} {\bf 85} (1967), 218-245.

\vskip .1in

\noindent [32] Ralph M. Kaufmann and R. C. Penner, ``Closed/open string diagrammatics'', Preprint math.GT/0603485, to appear {\it Nucl. Phys. B}.

\bye